\documentclass[a4paper]{amsart} 
\usepackage{amssymb} 
\usepackage{verbatim}
\usepackage[all]{xy}
\usepackage{multirow}
\usepackage{enumerate}

\calclayout

\title{Cofinite Submodule Closed Categories and the Weyl group}

\author{Apolonia Gottwald}

\address{Apolonia Gottwald, Fakult\"at f\"ur Mathematik, Universit\"at Bielefeld, D-33501 Bielefeld, Germany}

\email{alogisma@gmx.de}

\newtheorem{l1}{Lemma}[section]
\newtheorem{c1}[l1]{Corollary}
\newtheorem{t1}[l1]{Theorem}
\newtheorem{p1}[l1]{Proposition}

\theoremstyle{definition}
\newtheorem{d1}[l1]{Definition}
\newtheorem{e1}[l1]{Example}

\theoremstyle{remark}
\newtheorem{r1}[l1]{Remark}

\renewcommand{\mod}{\operatorname{\mathsf{mod}}\nolimits}
\newcommand{\ind}{\operatorname{\mathsf{ind}}\nolimits}

\newcommand{\Ker}{\operatorname{\mathsf{Ker}}\nolimits}
\newcommand{\Ext}{\operatorname{\mathsf{Ext}}\nolimits}
\newcommand{\End}{\operatorname{\mathsf{End}}\nolimits}
\renewcommand{\Im}{\operatorname{\mathsf{Im}}\nolimits}

\newcommand{\Hom}{\operatorname{\mathsf{Hom}}\nolimits}

\begin{document}
	
	\begin{abstract}
		We consider hereditary Artin algebras over arbitrary fields and prove that there is a natural bijection between the Weyl groups and the sets of full additive cofinite submodule closed subcategories of the module categories. While Oppermann, Reiten and Thomas have shown this for algebraically closed fields and finite fields, we give a different method of proof that holds independently of the field. 
		
		In particular, we show a relatively simple way to construct all modules that contain a given preinjective module as a submodule.
	\end{abstract}
	
	\date{\today}
	\maketitle
	
	\tableofcontents	
	
	\section{Introduction} \label{Intro}
	
	While submodule closed subcategories have not yet been extensively studied, there are many connections to different parts of representation theory. For example, if $A$ is a finite dimensional algebra, then every infinite submodule closed subcategory of $\mod A$ contains a minimal infinite submodule closed category, see \cite{Rin}.
	
	Submodule closed subcategories can also be used to prove that there is a filtration of the Ziegler spectrum that is indexed by the Gabriel-Roiter filtration, see \cite{KP}.

	In this paper, let $A$ be a hereditary Artin algebra over an arbitrary field.
		
	We aim to prove that there is a natural bijection between the Weyl group and the set of full additive cofinite submodule closed subcategories of the module category. Oppermann, Reiten and Thomas have shown this in \cite{ORT} for algebraically closed fields and finite fields. While we use the same bijection, we will give a completely different method of proof that does not depend on the field.
	
	Our paper is organized in the following way:
	
	First of all, we regard the Weyl group as a Coxeter group, see Section~\ref{SecCoxeter}. This allows us to regard the Weyl group elements as equivalence classes of words. In Section~\ref{SecWeylGroup}, we define a total order on these words and call the smallest word of each equivalence class \textit{leftmost}. Then we collect some results about this order.
	
	We conclude Section~\ref{SecWeylGroup} by stating the bijection, which is induced by a map between words of Weyl group elements and sets of preinjective modules. In Section~\ref{SecProof} to Section~\ref{SecProof3}, we will prove that a cofinite, full additive subcategory is submodule closed if and only if a leftmost word is mapped to its complement. Since we can assign a unique leftmost word to every element of the Weyl group, this gives a bijection between the full additive cofinite submodule closed subcategories and the Weyl group.
	
	For this proof, we will use the results of Section~\ref{SecPreinjModules}, which is devoted to monomorphisms between preinjective modules. In particular, we give a way to construct all modules that contain a given preinjective module as a submodule. This allows us to draw some lemmas in Section~\ref{SecConnection} about the structure of full additive cofinite submodule closed subcategories and how they are related to the words of Weyl group elements.
	
	In the sections \ref{SecProof} to \ref{SecProof3}, we use this to prove inductively that the proposed bijection exists. Finally, we conclude this chapter with some corollaries.
	A submodule closed subcategory will in the following always denote a full additive submodule closed subcategory of $\mod A$.
	
	\section{The Weyl group as a Coxeter group} \label{SecCoxeter}
	
	We define words following \cite{Lot}, pp. 1-4: 
	\begin{d1}
		Let $S$ be a set. We call $S$ an \textit{alphabet} and its elements \textit{letters}. A \textit{word} over the alphabet $S$ is a finite sequence 
		\begin{equation*}
		(s_1, s_2, \dots, s_n), s_i \in S.
		\end{equation*}
		The product of two words 
		is just the concatenation of the sequences. This product is associative and by identifying a letter $s \in S$ with the sequence $(s)$, we can write the word $(s_1, s_2, \dots, s_n)$ as the product $s_1s_2 \dots s_n$. The neutral element for this product is the empty word, which we accordingly denote as 1. 
		Thus, the set of words over $S$ together with the concatenation forms a \textit{monoid} $S^*$.
		
		If $\underline{w}:= s_1s_2 \dots s_n$ is a word over $S$, then $l(\underline{w}) := n$ is called the \textit{length} of $\underline{w}$. Furthermore, a word of the form $\underline{v} = s_{i_1}s_{i_2} \dots s_{i_m}$ with 
		\begin{equation*}
		1 \le i_1 < i_2 < \dots < i_m \le n
		\end{equation*}
		and $m \le n$ is a \textit{subword} of $\underline{w}$.
		
		If $\underline{v} = s_1 s_2 \dots s_m$ with $m \le n$, then we say that $\underline{v}$ is an \textit{initial subword} of $\underline{w}$.
	\end{d1}
	
	An introduction into Coxeter groups can be found in \cite{BB}. First, we need the definition of a Coxeter group, see pp. 1-2:
	
	\begin{d1}
		Let $S$ be a set and $W$ a group generated by $S$. Then $W$ is called a \textit{Coxeter group} if all relations have the form $(ss')^{m(s,s')} = 1$ with $s, s' \in S$ so that
		\begin{enumerate}
			\item $m(s,s') = 1$ if and only if $s=s'$.
			
			\item If $m(s,s')$ exists, then $m(s', s)$ also exists and $m(s,s') = m(s',s)$.
		\end{enumerate}
		If there is no relation between $s$ and $s'$, then we write $m(s,s') = m(s',s) = \infty$.
	\end{d1}
	
	We can describe the Coxeter group $W$ through the monoid $S^*$, see \cite{BB}, p. 3:
	
	\begin{p1} \label{LemWeylWords}
		Let $S$ be a set and $S^*$ the monoid of words over $S$. Let $W$ be a Coxeter group generated by $S$ with relations $(ss')^{m(s,s')}=1$. 
		
		Set $\equiv$ to be the equivalence relation on $S^*$ which is generated by allowing the insertion or deletion of words of the form
		\begin{equation*}
		(ss')^{m(s,s')} = \underbrace{ss'ss' \dots ss'}_{2m(s,s') \text{ letters}} 
		\end{equation*}
		for all $m(s,s') < \infty$.
		Then $S^*/\equiv$ is isomorphic to $W$.
	\end{p1}
	
	We will use the following notation:
	
	\begin{d1}
		Set $\{ss'\}^a := \underbrace{ss'ss' \dots}_{a \text{ letters}}$.
	\end{d1}
	
	The next lemma makes it easier to work with the relations, see \cite{BB}, p. 2:
	\begin{l1} 
		Let $S$, $W$ and $\equiv$ be as in Proposition \ref{LemWeylWords} and $s, s' \in S$. The equivalence of words $\{s's\}^{a} \equiv \{ss'\}^{a}$ holds if and only if $m(s,s')$ is a factor of $a$.
	\end{l1}

	Let $S_1, \dots S_n$ with $n \in \mathbb{N}$ be a complete list of non-isomorphic simple modules of the Artin algebra $A$.
	
	We can associate to $A$ a Cartan matrix, as in \cite{ARS}, pp 69, 241 and 288:
	
	\begin{d1} \label{DefCartan}
		To a hereditary Artin algebra $A$ we associate the \textit{Cartan matrix} $C = (c_{ij})_{nn}$ of the underlying graph of the quiver $A^{op}$.
		
		That is, we set $c_{ii} = 2$. If $i \neq j$ and $\Ext^1(S_i, S_j) = \Ext^1(S_j, S_i) = 0$, then $c_{ij} = c_{ji} = 0$. Finally, if $\Ext^1(S_i, S_j) \neq 0$, set 
		\begin{equation*}
		c_{ij} = - \dim_{\End_A (S_i)^{op}} \Ext^1(S_i, S_j)
		\end{equation*} 
		and 
		\begin{equation*}
		c_{ji} = - \dim_{\End_A (S_j)} \Ext^1(S_i, S_j).
		\end{equation*}
	\end{d1}
	
	A description of the Weyl group as a Coxeter group can be found in \cite{Kac}, Proposition 3.13:
	
	\begin{p1} \label{LemCoxeter}
		The Weyl group associated to $A$ with the Cartan matrix $(c_{ij})_{nn}$ is a Coxeter group generated by the reflections $s_1, s_2, \dots, s_n$ with relations $s^2_i = 1$ for all $1 \le i \le n$ and $(s_i s_j)^{m_{ij}} = 1$ for all $i \neq j$, where $m_{ij}$ depends on $c_{ij}c_{ji}$ in the following way:
		
		\begin{equation*}
		\begin{array}{l|lllll} c_{ij}c_{ji} & 0 & 1 & 2 & 3 & \ge 4\\ \hline
		m_{ij} & 2 & 3 & 4 & 6 & \infty
		\end{array}.
		\end{equation*}
		
	\end{p1}
	
	We can write all relations as $(s_is_j)^{m_{ij}}$ if we set $m_{ij} := 1$ for $i = j$. 
	
	Every element of the Weyl group is the equivalence class of several different words over the alphabet $S := \{s_1, s_2, \dots, s_n\}$. To distinguish between the elements of the Weyl group and the words over $S$, we will always use underlined letters to denote words and normal letters for Weyl group elements.
	
	\begin{r1}
		For $A = kQ$ with a field $k$ and a quiver $Q$ without oriented cycles, the relations depend only on the edges in the underlying graph of $Q$, see e.g. \cite{ORT}, p. 570:
		
		We have $m_{ij} = 2$ if there is no edge between the vertices $i$ and $j$ and $m_{ij} = 3$ if there is exactly one edge between $i$ and $j$. If there are two or more edges between $i$ and $j$, then $m_{ij} = \infty$.
	\end{r1}
	
	\begin{e1} \label{ExWeyl}
		Let $Q$ be the quiver
		\[
		\xymatrix @!0{1 \ar[r] \ar[d] & 3 \\
			4 & 2 \ar[l] \ar[u] }.
		\]
		The Weyl group of $A = kQ$ is a Coxeter group with the relations $(s_i s_j)^{m_{ij}} = 1$ and the following values for $m_{ij}$:
		
		\begin{tabular}{l|llll}  & 1 & 2 & 3 & 4 \\ \hline
			1 & 1 & 2 & 3 & 3 \\
			2 & 2 & 1 & 3 & 3 \\
			3 & 3 & 3 & 1 & 2 \\
			4 & 3 & 3 & 2 & 1
		\end{tabular}.
	\end{e1}	
	
	\section{Leftmost words} \label{SecWeylGroup}
	
	Let $A$ be a hereditary Artin algebra and $\mod A$ the category of finitely presented modules over $A$. Furthermore, let $\mathcal{I}$ be the subcategory of $\mod A$ consisting of all preinjective modules.
	
	Let $W$ be the Weyl group of $A$ and $S := \{s_1, s_2, \dots, s_n\}$ be the set of generators of $W$. Furthermore, let the simple modules $S_1, \dots S_n$ of $A$ with injective envelopes $I_1, \dots I_n$ be ordered in such a way that $\Hom(I_i, I_j) = 0$ if $i < j$. This is possible because $A$ is a hereditary Artin algebra, see \cite{ARS}, Chapter VIII, Proposition 1.5.
	
	\begin{d1} \label{DefWord}
		Consider $\mathcal{N} = (\mathbb{N}_0 \times \{1,2, \dots, n\}, <)$, where $<$ is the lexicographic order: for pairs $(r, i), (r', j) \in \mathcal{N}$, we have $(r, i) < (r', j)$
		if and only if one of the following holds:
		\begin{enumerate}
			\item $r < r'$
			
			\item $r = r'$ and $i < j$.
		\end{enumerate}
		Let $\underline{w} = s_{i_1} s_{i_2} \dots s_{i_m}$ be a word over the alphabet $S$ and $0 = r_1 \le  r_2 \le  \dots \le r_m \in \mathbb{N}_0$ the smallest non-negative integers so that 
		\begin{equation*}
		(r_1, i_1) < (r_2, i_2) < \dots < (r_m, i_m)
		\end{equation*}
		is fulfilled. Then we define 
		\begin{equation*}
		\rho(\underline{w}) := (r_1, i_1) (r_2, i_2) \dots (r_m, i_m).
		\end{equation*}
	\end{d1}
	
	\begin{e1}
		Consider the Weyl group of the quiver $Q$ from Example \ref{ExWeyl}. If we set $\underline{w}:=s_2s_3s_1s_3s_4s_1$ then
		\begin{equation*}
		\rho(\underline{w})= (0,2)(0,3)(1,1)(1,3)(1,4)(2,1).
		\end{equation*} 
	\end{e1}
	
	Now we can define a total order $<_l$ on the words of $W$; this is again a lexicographic order:
	
	\begin{d1} \label{DefOrder}
		Consider two words $\underline{w}, \underline{w}'$ with
		\begin{equation*}
		\rho(\underline{w}) = (r_1, i_1) (r_2, i_2) \dots (r_m, i_m)
		\end{equation*} 
		and 
		\begin{equation*}
		\rho(\underline{w}') = (r'_1, i'_1) (r'_2, i'_2) \dots (r'_{m'}, i'_{m'}).
		\end{equation*} 
		
		We write $\underline{w} <_l \underline{w}'$ if one of the following holds: 
		\begin{enumerate}
			\item $m < m'$
			\item $m = m'$ and there is a $j \in \mathbb{N}$ so that 
			\begin{equation*}
			(r_1, i_1) = (r'_1, i'_1), (r_2, i_2) = (r'_2, i'_2), \dots, (r_{j-1},i_{j-1}) = (r'_{j-1},i'_{j-1})
			\end{equation*}
			and
			\begin{equation*}
			(r_j,i_j) < (r'_j,i'_j).
			\end{equation*}
		\end{enumerate}
	\end{d1} 
	
	Now we define the leftmost word; this definition can be found for example in \cite{Arm}, p. 411:

	\begin{d1} \label{DefLeftmost}
		We call a word $\underline{w}$ for $w \in W$ \textit{leftmost} if for every other word $\underline{w}'$ for $w$ the inequality $\underline{w} <_l \underline{w}'$ holds.
	\end{d1}

	\begin{e1}
		For the Weyl group from Example \ref{ExWeyl}, the words
		\[
		s_3 <_l s_2 s_3 <_l s_3 s_2 <_l s_2s_3s_2
		\]
		are all leftmost words and
		\[
		s_2 s_3 s_2 <_l s_3 s_2 s_3 <_l s_2 s_3 s_1 s_2 s_1
		\]
		are all words for the same element of the Weyl group.
	\end{e1}
	
	Since $<_l$ is a total order, every element $w \in W$ has a unique leftmost word. Obviously, the leftmost word is \textit{reduced}, that is, it has the smallest possible length for a word of $w$.
	
	We follow with a Lemma about the order $<_l$ and the relations:
	
	\begin{l1} \label{LemOneExchange}
		Suppose that $\underline{w}_1 = \underline{u} \{s_is_j\}^{m_{ij}} \underline{v}$ 	for some reflections $s_i$, $s_j$, $i \neq j$, words $\underline{u}$, $\underline{v}$. Set
		\begin{equation*}
		\rho(\underline{w}_1) = \rho (\underline{u}) \underbrace{(p,i) (q,j) (p+1,i) \dots}_{m_{ij} \text{ pairs}} \rho_1
		\end{equation*}
		for some $p,q \in \mathbb{N}_0$ and a sequence of pairs $\rho_1$. Set 
		\begin{equation*}
		\underline{w}_2 = \underline{u} \{s_js_i\}^{m_{ij}} \underline{v}.
		\end{equation*}
		Then $\underline{w}_2 <_l \underline{w}_1$ if and only if both of the following conditions are fulfilled:
		\begin{enumerate}
			\item 
			$1 \le q$.
			
			\item Let $(r,k)$ be a pair in $\rho(\underline{u})$. Then $(r,k) < (q-1, j)$.
		\end{enumerate}
	\end{l1}
	
	\begin{proof}
		Let  $\rho(\underline{w}_2) = \rho(\underline{u})\underbrace{(q',j) (p',i) \dots}_{m_{ij} \text{ pairs}} \rho_2$ for some sequence of pairs $\rho_2$.
		
		Suppose that $\underline{w}_2 <_l \underline{w}_1$. Then $(q',j) < (p,i)$ by Definition \ref{DefOrder} and thus $q = q'+1$. So the first condition is fulfilled.
		
		Now consider a pair $(r,k)$ in $\rho(\underline{u})$. Then $(r,k) < (q',j) = (q-1,j)$ and the second condition is fulfilled.
		
		On the other hand, suppose that the conditions 1 and 2 are fulfilled. Then $q'$ is the smallest integer so that $(q',j)$ is bigger than all $(r,k)$ in $\rho(\underline{u})$. By the second condition, $(q',j) \le (q-1,j)$. 
		
		Furthermore, $q$ is the smallest integer so that $(p,i) < (q,j)$. It follows that $(q-1,j) < (p,i)$, since $i \neq j$.
		
		Together, $(q',j) \le (q-1,j)< (p,i)$ and by Definition \ref{DefOrder}, we have $\underline{w}_2 <_l \underline{w}_1$.
	\end{proof}
	
	The following lemmas are important for the induction with which we prove the main theorem of this chapter:

	\begin{l1} \label{LemSame}
		Let $\underline{x}, \underline{x}', \underline{y}$ be words and $s_i, s_j$ reflections. We suppose that the words $\underline{w}=\underline{x}s_i\underline{y}$ and $\underline{w}'= \underline{x}'s_j\underline{y}$ are equivalent, $\underline{x}$ is leftmost and $\underline{w} <_l \underline{w}'$. Let $\underline{z}$ be the longest initial subword that $\underline{w}$ and $\underline{w}'$ share. If there is no $\underline{x}'' \equiv \underline{x}'$ that shares an initial subword with $\underline{w}$ which is longer than $\underline{z}$, then there are pairs $(r, h), (s, i), (t,j)$ and series of pairs $\rho_1, \rho_2, \rho_3$, $\rho_4$ so that $\rho(\underline{x}') = \rho_1(r, h)\rho_2$
		\begin{equation*}
		\rho(\underline{w}) = \rho_1 (r, h) \rho_2 (s, i) \rho_3
		\end{equation*}
		and there is some word $\underline{w}'' \equiv \underline{w}$ with
		\begin{equation} \label{Form7}
		\rho(\underline{w}'') = \rho_1 \rho_2 (s, i) (t, j)\rho_4.
		\end{equation}
		Either $\rho_3 = \rho_4$, or a pair $(q, g)$ is in $\rho_4$ if and only if $(q-1, g)$ is in $\rho_3$.
		
		If $\underline{v}$ is the initial subword of $\underline{w}$ with $\rho(\underline{v})=\rho_1$, then no relation on reflections in $\underline{v}$ is needed to transform $\underline{w}$ into $\underline{w}''$.
	\end{l1}
	
	\begin{proof}
		We prove this inductively and assume without loss of generality that $m_{ij}$ is odd. If $m_{ij}$ is even, we only need to relabel $s_i$ and $s_j$ in the arguments below.
		
		If there is some word $\underline{x}_1$ so that 
		\begin{equation} \label{WForm}
		\underline{w} = \underline{x}_1\{s_is_j\}^{m_{ij}}\underline{y} <_l \underline{x}_1\{s_js_i\}^{m_{ij}}\underline{y}
		\end{equation}
		then $\underline{w}'' :=  \underline{x}_1\{s_js_i\}^{m_{ij}}\underline{y}$ fulfils the assertions by Lemma \ref{LemOneExchange}. 
		
		On the other hand, we get a similar result if 
		\begin{equation*}
		\underline{x}_1\{s_js_i\}^{m_{ij}}\underline{y} <_l \underline{x}_1\{s_is_j\}^{m_{ij}}\underline{y}:
		\end{equation*}
		Then
		\begin{equation*}
		\rho(\underline{x}_1\{s_js_i\}^{m_{ij}}) = \rho_1 \rho'_2 (s-1, i) (t-1, j)\rho'_4,
		\end{equation*}
		where a pair $(q-1, g)$ is in $\rho'_2$ if and only if $(q, g)$ is in $\rho'_2$. Furthermore, either $\rho'_4 = \rho_3$ or a pair $(q-1, g)$ is in $\rho'_4$ if and only if $(q, g)$ is in $\rho_3$.
		
		Now suppose that $\underline{w} = \underline{u}\{s_ls_k\}^{m_{kl}}\underline{v}$ is not of the form in equation (\ref{WForm}), but the assertion is true for $\underline{w}_1 = \underline{u}\{s_ks_l\}^{m_{kl}}\underline{v}$.	
		
		If $\underline{w}_1 <_l \underline{w}'$, then there is some word $\underline{w}''_1$ with 
		\begin{equation*}
		\rho(\underline{w}_1) = \rho'_1 (r', h') \rho'_2 (s, i) \rho'_3
		\end{equation*}
		and 
		\begin{equation*}
		\rho(\underline{w}''_1) = \rho'_1 \rho'_2 (s, i) (t, j)\rho'_4.
		\end{equation*}
		for some pair $(r', h')$ and series of pairs $\rho'_1, \rho'_2, \rho'_3, \rho'_4$. Either $\rho'_3 = \rho'_4$, or a pair $(q, g)$ is in $\rho'_4$ if and only if $(q-1, g)$ is in $\rho'_3$.
		
		There is some $\underline{u}'$ so that $\rho'_1 = \rho(\underline{u}')$. 
		
		Furthermore, there is a pair $(q, l)$ and a series of pairs $\rho$ so that 
		\begin{equation*}
		\rho(\underline{w}) = \rho(\underline{u}) (q, l) \rho.
		\end{equation*}

		If $\underline{u}' = \underline{u}$ or $\underline{u}' = \underline{u} \{s_l s_k\}^{m_{lk}-1}$, then we set $(r, h) := (q, l)$. Since $\underline{x}$ is leftmost , $\underline{w} <_l \underline{w}_1$ and by Lemma \ref{LemOneExchange}, the assertion is true.
		
		Since $\underline{x}$ is reduced, there is only one other case: the word $\underline{w}''_1$ has $\{s_ks_l\}^{m_{kl}}$ as a subword and the relation $\{s_ks_l\}^{m_{kl}} \equiv \{s_ls_k\}^{m_{kl}}$ gives a word $\underline{w}''$ that fulfils (\ref{Form7}).
		
		It remains to prove the assertion if $\underline{w}' <_l \underline{w}_1$. 
		
		Then we can inductively assume that there is some $\underline{w}''_1$ so that
		\begin{equation*}
		\rho(\underline{w}''_1) = \rho_1 \rho'_2 (s-1, i) (t-1, j)\rho'_4,
		\end{equation*}
		where a pair $(q-1, g)$ is in $\rho'_2$ if and only if $(q, g)$ is in $\rho'_2$. Furthermore, either $\rho'_4 = \rho_3$ or a pair $(q-1, g)$ is in $\rho'_4$ if and only if $(q, g)$ is in $\rho_3$.
		
		But there is no $\underline{x}'' \equiv \underline{x}'$ that shares an initial subword with $\underline{w}$ which is longer than $\underline{z}$, the longest initial subword that $\underline{w}$ and $\underline{w}'$ share. 
		
		Since $\underline{w}''_1$ and $\underline{w}_1$ share the initial subword $\underline{u}'$ with $\rho(\underline{u}_1) = \rho_1$, this is only possible if $\underline{u}' = \underline{u} \{s_l s_k\}^{m_{lk}-1}$.
		
		Again, we set $(r, h) := (q, l)$. Since $\underline{x}$ is leftmost , $\underline{w} <_l \underline{w}_1$ and by Lemma \ref{LemOneExchange}, the assertion is true.
	\end{proof}
	
	Completely analogously, we can prove the following:
	
	\begin{l1} \label{LemSame2}
		Let $\underline{x}, \underline{x}', \underline{y}$ be words and $s_i, s_j$ reflections. If the words $\underline{w}=\underline{x}s_i\underline{y}$ and $\underline{w}'= \underline{x}'s_j\underline{y}$ are equivalent, $\underline{x}$ is leftmost and $\underline{w}' <_l \underline{w}$, then there are pairs $(r, h), (s, i)$ and series of pairs $\rho_1, \rho_2, \rho_3$, $\rho_4$ so that $\rho(\underline{x}') = \rho_1\rho_2$
		\begin{equation*}
		\rho(\underline{w}) = \rho_1 \rho_2 (s, i) \rho_3
		\end{equation*}
		and there is some word $\underline{w}'' \equiv \underline{w}$ with
		\begin{equation*}
		\rho(\underline{w}'') = \rho_1 (r, h) \rho_2 \rho_4.
		\end{equation*}
		Either $\rho_3 = \rho_4$, or a pair $(q, g)$ is in $\rho_4$ if and only if $(q+1, g)$ is in $\rho_3$.
		
		If $\underline{v}$ is the initial subword of $\underline{w}$ with $\rho(\underline{v})=\rho_1$, then no relation on reflections in $\underline{v}$ is needed to transform $\underline{w}$ into $\underline{w}''$.
	\end{l1}
	
	We get the following corollary:
	
	\begin{c1}   \label{LemIndStep}
		If $\underline{u}\{s_is_j\}^{m_{ij}-1}$ and $\underline{u}s_j$ are leftmost, then either $\underline{u}\{s_is_j\}^{m_{ij}}$ is leftmost or $\underline{u}s_j < \underline{u}s_i$. 
	\end{c1}
	
	\begin{proof}
		If we have  $\underline{u}\{s_js_i\}^{m_{ij}} <_l\underline{u}\{s_is_j\}^{m_{ij}}$, then $\underline{u}s_j < \underline{u} s_i$.
		
		Furthermore, if $\underline{u}\{s_is_j\}^{m_{ij}}$ is not leftmost, but $\underline{u}\{s_is_j\}^{m_{ij}-1}$ is, then we can write $\underline{w}=\underline{x}s_i\underline{y}$ and there is some $\underline{w}'= \underline{x}'s_j\underline{y}$ so that $\underline{w} \equiv \underline{w}'$, $\underline{x}$ is leftmost and $\underline{w}' <_l \underline{w}$.
		
		By Lemma \ref{LemSame2}, there are some words $\underline{u}_1, \underline{u}_2$ and a reflection $s_h$ so that
		\begin{equation*}
		\underline{u} = \underline{u}_1 \underline{u}_2
		\end{equation*}
		\begin{equation*}
		\underline{w}'' = \underline{u}_1 s_h \underline{u}_2 \{s_is_j\}^{m_{ij}-1} \equiv \underline{w}
		\end{equation*}
		and $\underline{w}'' <_l \underline{w}$.
		
		So $\underline{u}s_j \equiv \underline{u}_1 s_h \underline{u}_2$ and $\underline{u}_1 s_h \underline{u}_2<_l \underline{u}s_j$. 
	\end{proof}

	\begin{r1} \label{RemIndStep}
		Note that in Lemma \ref{LemSame2}, we do not actually need to assume that $\underline{x}$ is leftmost; it is sufficient that $\underline{x}$ is reduced and the following holds: let $\underline{x}'' \equiv \underline{x}$ so that $\underline{x}'' <_l \underline{x}$. Furthermore, assume that $\underline{x}_2$ is the maximal initial subword that $\underline{x}'$ and $\underline{x}$ share. Then there is some $\underline{w}' \equiv \underline{w}$ with $\underline{w}' <_l \underline{x}'s_i \underline{y}$.
		
		So analogously to \ref{LemIndStep}, we see: Suppose that there is some word $\underline{w}'$ so that
		$\underline{u}\{s_is_j\}^{m_{ij}} \equiv \underline{w}'$ and for all $\underline{u}' \equiv \underline{u}\{s_is_j\}^{m_{ij}-1}$, we have $\underline{w}' <_l \underline{u}'s_i$ if $m_{ij} = 3$ and $\underline{w}' <_l \underline{u}'s_j$ otherwise. Then $\underline{u}s_j$ is not leftmost.
	\end{r1}

	\begin{l1} \label{LemPart1}
		Suppose that $\underline{w} \equiv \underline{u}\{s_is_j\}^{m_{ij}}$ with $m_{ij} \ge 3$. If there is some $i \neq j \neq k$ with $\underline{w} = \underline{u}'s_k\{s_is_j\}^{m_{ij} - m}$ for some even $m \ge 2$ or $\underline{w} = \underline{u}'s_k\{s_js_i\}^{m_{ij} - m}$ for some odd $m \ge 3$, then $m_{ik} = 2$ or $m_{jk} = 2$.
	\end{l1}
	
	\begin{proof}
		This is a simple, inductive proof: without loss of generality, we can assume that $\underline{w} = \underline{u}_1s_k\{s_is_j\}^{m_{ij}-2}$ and $m_{jk} = 3$. Then there is some $\underline{u}_2$ so that
		\begin{equation*}
		\underline{w} \equiv\underline{u}_2s_ks_js_k\{s_is_j\}^{m_{ij}-2} \equiv \underline{u}_2s_js_k\{s_js_i\}^{m_{ij}-1}.
		\end{equation*}
		So there is some $\underline{u}_3$ so that
		\begin{equation*}
		\underline{u}_2s_js_k \equiv \underline{u}_3 \{s_is_k\}^{m_{ik}}.
		\end{equation*}
		If $m_{ik} \ge 3$, then we have the same situation as before, only considering a shorter word. Since $\underline{w}$ is finite, we see with induction on the length of $\underline{w}$ that $m_{ik} = 2$ or $m_{jk} = 2$.
	\end{proof}
	
	Similarly, we can prove the following:
	
	\begin{l1} \label{LemPart2}
		Suppose that $\underline{u}s_j$ is leftmost, but $\underline{u}s_js_i$ is not leftmost for some word $\underline{u}$ and some reflections $s_i, s_j$ with $m_{ij} \ge 4$.
		
		Then there are $s, t \in \mathbb{N}$ so that $\rho(\underline{u}s_js_i) = \rho(\underline{u})(t, j)(s, i)$ and $\rho(\underline{u})$ contains the pair $(s-1, i)$. If $m_{ij} = 6$, then $\rho(\underline{u})$ additionally contains the pairs $(s-2, i)$ and $(t-1, j)$.
	\end{l1}
	
	\begin{proof}
		Suppose that the assertions are not fulfilled.
		
		We can without loss of generality assume that there is some $\underline{u}'$ and reflections $s_{k_1}, \dots, s_{k_m}, s_{l_1}, \dots, s_{l_{m'}}$ so that
		\begin{equation*}
		\underline{u}s_js_i = \underline{u}'s_js_{k_1} \dots s_{k_m}s_i s_{l_1} \dots s_{l_{m'}}\{s_js_i\}^{m_{ij}-2} =: \underline{u}''\{s_js_i\}^{m_{ij}-2} 
		\end{equation*}
		Then there are $s', t', t'', q_1, \dots, q_m, r_1, \dots, r_{m'} \in \mathbb{N}$ so that
		\begin{equation*}
		\rho(\underline{u}''s_j) = \rho(\underline{u}')(t', j)(q_1, k_1) \dots (q_m, k_m)(s', i) (r_1, l_1) \dots (r_{m'}, l_{m'})(t'', j).
		\end{equation*}
		If $(t''-1, j) < (s', i)$, then the assertions of the lemma are fulfilled. Otherwise, we get one of the following cases:
		\begin{enumerate}[(a)]
			\item 	There is some $1 \le o \le m'$ with $(t''-1, j) < (r_o, l_o)$ with $m_{l_o, j} \ge 3$.
			
			\item The words $s_i s_{l_1} \dots s_{l_{m'}}s_j$ and $\underline{u}s_i$ are not leftmost, contrary to the assumptions.
		\end{enumerate}
		So we can assume that the first case is fulfilled. Furthermore, without loss of generality, we can assume $o = m'$.
		
		If $m_{il_1} = \dots  = m_{im'} = 2$, then
		\begin{equation*}
		\underline{u}'' \equiv \underline{u}' s_js_{k_1} \dots s_{k_m}s_{l_1} \dots s_{l_{m'}}s_i =: \underline{u}'''.
		\end{equation*}
		Either $\underline{u}''' <_l \underline{u}''$, contrary to the assumptions, or there is some $s'' > s'$ so that
		\begin{equation*}
		\rho(\underline{u}''') = \rho(\underline{u}')(t', j)(q_1, k_1) \dots (q_m, k_m)(r_1, l_1) \dots (r_{m'}, l_{m'})(s'', i).
		\end{equation*}
		Since $\underline{u}s_j$ is leftmost, but $\underline{u}s_js_i$ is not, there is some $v_1$ so that
		\begin{equation} \label{EqRel}
		\underline{u}s_js_i \equiv \underline{v}_1\{s_js_i\}^{m_{ij}}.
		\end{equation} 
		Thus, 
		\begin{equation*}
		\underline{u}'' \equiv \underline{u}' \{s_j s_{l_{m'}}\}^{m_{jl_{m'}}}s_i
		\end{equation*}
		and we are in an analogous situation to before, only considering a shorter word. Since the length of $\underline{w}$ is finite, the assertions of the lemma are inductively true under these assumptions.
		
		So we can assume without loss of generality that $m_{il_1} \ge 3$. Furthermore, we have $(s'+1, i)< (t''-1, j) < (r_{m'}, l_{m'})$. So there is at least one $1 \le o \le m'$ so that $m_{ol_{m'}} \ge 3$, since otherwise $s_js_{k_1} \dots s_{k_m}s_i s_{l_1}$ is equivalent to a smaller word (that begins with $s_j$), contrary to the assumption that $\underline{u}s_j$ is leftmost.
		
		Then there is some $\underline{v}_2$ so that
		\begin{equation*}
		\underline{u}'s_js_{k_1} \dots s_{k_m}s_i s_{l_1} \equiv \underline{v}_2\{s_is_{l_1}\}^{m_{il_1}}.
		\end{equation*}
		
		Because of (\ref{EqRel}), Lemma \ref{LemSame}, Lemma \ref{LemSame2} and $m_{il_1} \ge 3$, we get some $\underline{v}_3$ so that
		\begin{equation*}
		\underline{u}'' \equiv \underline{v}_3 s_js_is_{l_1} \dots s_{l_m}s_i.
		\end{equation*}
		Since $m_{jl_m} \ge 3$ and $m_{l_ol_m} \ge 3$, we are in the same situation as before, only considering a word of shorter length. Inductively, the proof is complete. 
	\end{proof}
	
	Now we can define an assignment which maps the words of the Weyl group to the cofinite full additive subcategories of $\mod A$. We will show that this map yields a bijection between the Weyl group and the set of cofinite submodule closed subcategories.
	
	Let $\tau = D \text{Tr}$ be the Auslander-Reiten translation, see \cite{ARS}, p.~106. With \cite{ARS}, p.~259, every indecomposable preinjective module is of the form $\tau^r I_i$ for some $r \in \mathbb{N}$ and $1 \le i \le n$.
	
	\begin{d1} \label{DefBijection}
		We can identify the pairs in $\mathcal{N}$ and the indecomposable preinjective modules by setting $(r,i) = \tau^r I_i$.
		
		Not only does this give us a natural order on the preinjective modules, but this also yields an injective map from the words of the Weyl group to the cofinite full additive subcategories of $\mod A$: If  
		\begin{equation*}
		\rho(\underline{w}) = (r_1, i_1) (r_2, i_2) \dots (r_m, i_m),
		\end{equation*} 
		then $\underline{w} \mapsto \mathcal{C}_{\underline{w}}$, where $\mathcal{C}_{\underline{w}}$ is the full additive category with
		\begin{equation*}
		\ind \mathcal{C}_{\underline{w}}= \ind A \setminus \{ (r_1, i_1), (r_2, i_2), \dots, (r_m, i_m)\}.
		\end{equation*}
		
		For a Weyl group element $w$ with the leftmost word $\underline{w}$, define $\mathcal{C}_{w} := \mathcal{C}_{\underline{w}}$.
	\end{d1}

	\begin{e1}
		Let $A$ be as in Example \ref{ExWeyl} and $\underline{w}= s_1s_2s_3s_2s_4s_1$. Then
		\[
		\rho(\underline{w}) = (0,1)(0,2)(0,3)(1,2)(1,4)(2,1)
		\] 
		and
		\begin{align*}
		\ind \mathcal{C}_{\underline{w}} = \ind A \setminus \{I_1, I_2, I_3, \tau I_2, \tau I_4, \tau^2 I_1\}.
		\end{align*}
	\end{e1}
	
	We will prove that the restriction of this map on the leftmost words is a bijection between those and the cofinite submodule closed subcategories.
	
	Since every element of the Weyl group has a unique leftmost word, this gives a bijection between the elements of the Weyl group and the cofinite submodule closed subcategories.
	
	The same bijection is used in \cite{ORT}.

	\section{Monomorphisms between preinjective modules} \label{SecPreinjModules}
	An observation makes the aim of the chapter much simpler to achieve: the cofinite submodule closed subcategories of the module category correspond naturally to the cofinite submodule closed subcategories of $\mathcal{I}$, the category of the preinjective modules. 
	
	Thus we devote this section to preinjective modules. In particular, we give a way to construct all modules $U$ that contain a given preinjective, indecomposable module $M$ as a submodule.
	
	In Section \ref{SecConnection} we will use this to show the connection to the Coxeter structure of the Weyl group. In Section \ref{SecProof} to \ref{SecProof3}, we will use this connection to prove that the bijection that we described exists.
	
	\begin{p1} \label{LemPreinj}
		There is a bijection between full additive cofinite submodule closed subcategories of $\mod A$ and full additive cofinite submodule closed subcategories of $\mathcal{I}$. It maps the category $\mathcal{C}$ to the category $\mathcal{C'} = \mathcal{C} \cap \mathcal{I}$. Furthermore, \begin{equation*}
		\ind A \setminus \mathcal{C} = \ind \mathcal{I} \setminus \mathcal{C}.
		\end{equation*}
	\end{p1}
	
	\begin{proof}
		This is completely analogous to \cite{ORT}, Proposition 2.2: 
		
		If $A$ is representation finite, then $\mod A = \mathcal{I}$ and there is nothing to prove. Suppose that $A$ is not representation finite. Since $\mathcal{C}$ is cofinite, there is some $r \in \mathbb{N}$, so that $\tau^r I_1, \tau^r I_2, \dots, \tau^r I_n \in \mathcal{C}$. Now suppose that $M$ is a preprojective or regular module. Then $\tau^{-r} M$ exists and has an injective envelope $I$. Since $\tau^r$ preserves monomorphisms, $M \subseteq \tau^r I \in \mathcal{C}$. So $M \in \mathcal{C}$ and $\ind A \setminus \mathcal{C} = \ind \mathcal{I} \setminus \mathcal{C}$.
		
		Thus the assignment $\mathcal{C} \mapsto \mathcal{C} \cap \mathcal{I}$ is a bijection between the full additive cofinite submodule closed subcategories of $\mod A$ and the full additive cofinite submodule closed subcategories of $\mathcal{I}$.
	\end{proof}
	
	We start the construction of exact sequences with a lemma that holds for all Artin algebras:

	\begin{l1} \label{LemARgen}
		Let $A$ be an arbitrary Artin algebra and $M, X \in \mod A$ indecomposable. Let
		\begin{equation} \label{ExSeq}
		\xymatrixcolsep{3pc}\xymatrix{0 \ar[r] & M \ar^-{\left[\begin{smallmatrix}
				f_1\\f_2
				\end{smallmatrix}\right]}[r] & X \oplus X' \ar^-{\left[\begin{smallmatrix}
				g_{11} & g_{12}\\0& g_{22}
				\end{smallmatrix}\right]}[r] & Y \oplus Y' \ar[r] & 0}
		\end{equation}
		be an exact sequence for some $X',Y,Y' \in \mod A$, so that there is some $Z \in \mod A$ and an AR-sequence
		\begin{equation} \label{ExSeq2}
		\xymatrixcolsep{3pc}\xymatrix{0 \ar[r] & X \ar^-{\left[\begin{smallmatrix}
				g_{11}\\f'_2
				\end{smallmatrix}\right]}[r] & Y \oplus Z \ar^-{\left[\begin{smallmatrix}
				g'_{1} & g'_{2}
				\end{smallmatrix}\right]}[r] & \tau^{-1} X \ar[r]& 0}.
		\end{equation}
		
		If for some $U \in \mod A$, a monomorphism $h: M \rightarrowtail U$ factors through $f = \left[\begin{smallmatrix}
		f_1\\f_2
		\end{smallmatrix}\right]$ and $X \nmid U$, then $h$ also factors through $f'' = \left[\begin{smallmatrix}
		-f'_2f_1\\f_2
		\end{smallmatrix}\right]$
		and the following sequence is exact:
		\begin{equation*}
		\xymatrixcolsep{4pc}\xymatrix{0 \ar[r] & M \ar^-{\left[\begin{smallmatrix}
				-f'_2f_1\\f_2
				\end{smallmatrix}\right]}[r] & Z \oplus X' \ar^-{\left[\begin{smallmatrix}
				g'_{2} & g'_1g_{12}\\0& g_{22}
				\end{smallmatrix}\right]}[r] & \tau^{-1} X \oplus Y' \ar[r] & 0}.
		\end{equation*}
	\end{l1}
	
	\begin{proof}
		By (\ref{ExSeq2}), the sequence
		\begin{equation*}
		\xymatrixcolsep{4pc}\xymatrix{0 \ar[r] & X \ar^-{\left[\begin{smallmatrix}
				g_{11}\\f'_2\\0
				\end{smallmatrix}\right]}[r] & Y \oplus Z \oplus Y' \ar^-{\left[\begin{smallmatrix}
				g'_{1} & g'_{2} & 0\\
				0&0&\text{id}_{Y'}
				\end{smallmatrix}\right]}[r] & \tau^{-1} X \oplus Y'\ar[r]& 0}
		\end{equation*}
		is also exact.  By \cite{ARS}, Chapter I, Corollary 5.7, the diagrams
		\begin{equation}  \label{Cartesian1}
		\xymatrix{M \ar^{f_1}[r]  \ar_{-f_2}[d]& X \ar^-{\left[\begin{smallmatrix}
				g_{11}\\0
				\end{smallmatrix}\right]}[d]\\
			X' \ar_-{\left[\begin{smallmatrix}
				g_{12}\\g_{22}
				\end{smallmatrix}\right]}[r] & Y \oplus Y'}
		\end{equation}
		and
		\begin{equation*}
		\xymatrixcolsep{3pc}\xymatrix{X \ar^{-f'_2}[r] \ar_-{\left[\begin{smallmatrix}
				g_{11}\\0
				\end{smallmatrix}\right]}[d]& Z \ar^-{\left[\begin{smallmatrix}
				g'_{2}\\0
				\end{smallmatrix}\right]}[d]\\
			Y \oplus Y' \ar_-{\left[\begin{smallmatrix}
				g'_{1}&0\\0 & \text{id}_{Y'}
				\end{smallmatrix}\right]}[r] & \tau^{-1} X \oplus Y'}
		\end{equation*}
		are both pushouts and pullbacks. So the diagram 
		\begin{equation*}
		\xymatrixcolsep{3pc}\xymatrix{M \ar^{-f'_2f_1}[r] \ar_-{-f_2}[d]& Z \ar^-{\left[\begin{smallmatrix}
				g'_{2}\\0
				\end{smallmatrix}\right]}[d]\\
			X' \ar_-{\left[\begin{smallmatrix}
				g'_{1}g_{12}\\g_{22}
				\end{smallmatrix}\right]}[r] & \tau^{-1} X \oplus Y'}
		\end{equation*}
		is itself a pushout and a pullback, see \cite{Gra}, p 334-335. By the definition of the Auslander-Reiten sequence (see \cite{ARS}, p. 137 and p. 144), the sequence
		\begin{equation*}
		\xymatrixcolsep{4pc}\xymatrix{0 \ar[r] & M \ar^-{\left[\begin{smallmatrix}
				-f'_2f_1\\f_2
				\end{smallmatrix}\right]}[r] & Z \oplus X' \ar^-{\left[\begin{smallmatrix}
				g'_{2} & g'_1g_{12}\\0& g_{22}
				\end{smallmatrix}\right]}[r] & \tau^{-1} X \oplus Y' \ar[r] & 0}
		\end{equation*}
		is exact. It remains to show that $h: M \rightarrowtail U$ factors through $f'' = \left[\begin{smallmatrix}
		-f'_2f_1\\f_2
		\end{smallmatrix}\right]$.
		
		Since we have assumed that $h$ factors through $f = \left[\begin{smallmatrix}
		f_1\\f_2
		\end{smallmatrix}\right]$, there is a morphism $s = \left[s_1 s_2\right]: X \oplus X' \rightarrow U$ so that 
		\begin{equation*}
		h = \left[\begin{smallmatrix}
		s_1&s_2
		\end{smallmatrix}\right] \left[\begin{smallmatrix}
		f_1\\f_2
		\end{smallmatrix}\right] = s_1f_1 +s_2f_2
		\end{equation*}
		By the definition of the Auslander-Reiten sequence (see \cite{ARS}, p. 137 and p. 144), the morphism $s_1: X \rightarrow U$ factors through $\left[\begin{smallmatrix}
		g_{11}\\f'_2
		\end{smallmatrix}\right]$: there is a morphism $s' = \left[s'_1 s'_2\right]: Y \oplus Z \rightarrow U$ so that
		\begin{equation*}
		s_1 = \left[\begin{smallmatrix}
		s'_1&s'_2
		\end{smallmatrix}\right] \left[\begin{smallmatrix}	g_{11}\\f'_2
		\end{smallmatrix}\right] = s'_1g_{11} +s'_2f'_2.
		\end{equation*}
		So we get
		\begin{equation*}
		h = s'_1g_{11}f_1 + s'_2f'_2f_1 + s_2f_2.
		\end{equation*}
		Since (\ref{Cartesian1}) is commutative, we have
		\begin{equation}
		h = -s'_1g_{12}f_2+ s'_2f'_2f_1 + s_2f_2 =  \left[\begin{smallmatrix}
		-s'_2 & s_2-s'_1g_{12}
		\end{smallmatrix}\right] \left[\begin{smallmatrix}	-f'_2f_1\\f_2
		\end{smallmatrix}\right]
		\end{equation}
		and $h$ factors through $f''$.
	\end{proof}

	We can even say more:
	
	\begin{l1} \label{LemExSeq}
		Let $A$ be a hereditary Artin algebra, $M \in \mathcal{I}$ and $U \in \mod A$. Suppose that the sequences of modules
		\begin{align*}
		&(X_1, X_2, \dots, X_m)\\
		&(X'_1, X'_2,\dots, X'_m)\\
		&(Y_1, Y_2, \dots, Y_m)
		\end{align*}
		fulfil the following conditions: 
		\begin{enumerate}[(S1)]
			\item There is an Auslander-Reiten sequence
			\begin{equation*}
			\xymatrix{0 \ar[r] & M \ar[r] & X_1 \oplus X'_1 \ar[r] & Y_1 \ar[r]& 0}.
			\end{equation*}
			
			\item For all $1 \le i < m$, there is some $\alpha_i \in \mathbb{N}$ so that $X_i^{\alpha_i} \mid X_i \oplus X'_i$, but $X_i^{\alpha_i} \nmid U$.
			
			\item For $1 \le i < m$, there is an Auslander-Reiten sequence of the form
			\begin{equation*}
			\xymatrix{0 \ar[r] & X_i \ar[r] & Z_i \ar[r] & \tau^{-1} X_i \ar[r]& 0}.
			\end{equation*}
			Let $Y'_i$ be the maximal module that is a direct summand of both $Y_i$ and $Z_i$. Write $Y_i = Y'_i \oplus Y''_i$ and $Z_i = Y'_i \oplus Z'_i$.
			
			If $\tau^{-1} X_i \mid X'_i$, then let $X''_i$ be the module so that $X'_i = \tau^{-1} X_i \oplus X''_i$ and set $Y'''_i := 0$. Otherwise, set $X''_i := X'_i$ and $Y'''_i := \tau^{-1} X_i$.

			The following equations hold:
			\begin{align*}
			X_{i+1} \oplus X'_{i+1} &= X''_i \oplus Z'_i\\
			Y_{i+1} &= Y''_i \oplus Y'''_i.
			\end{align*}
		\end{enumerate}
		Then for all $1 \le i \le m$ there is an exact sequence
		\begin{equation} \label{SeqAlg1}
		\xymatrix{0 \ar[r] & M \ar^-{f_i}[r] & X_i \oplus X'_i \ar^-{g_i}[r] & Y_i \ar[r]& 0}.
		\end{equation}
		Furthermore, if a monomorphism $M \rightarrowtail U$ exists, then it factors through all $f_i$. 
	\end{l1}
	
	To prove Lemma \ref{LemExSeq}, we need the following observation:
	
	\begin{r1} \label{RemObs}
		Suppose that 
		\begin{equation*}
		X_i \oplus X'_i = X_i \oplus X_{i+1} \oplus B_i
		\end{equation*} 
		and that $C_i$ is the maximal module that is both a direct summand of $Y_i$ and $Z_{i+1}$. Furthermore, write $Y_i = C_i \oplus C'_i$ and $Z_{i+1} = C_i \oplus D_{i}$. Set $B_i = B'_i$ and $C''_i = \tau X_i$ if $\tau X_i \nmid B_i$ and $B_i = B'_i \oplus \tau X_i$ and $C''_i = 0$ if $\tau X_i \mid B_i$. Then the sequences of modules 
		\begin{align*}
		&(X_1, X_2\dots, X_{i-1}, X_{i+1}, X_i, X_{i+2}, X_{i+3}, \dots, X_m)\\
		&(X'_1, X'_2\dots, X'_{i-1}, X_i \oplus B_i, B'_i\oplus D_i, X'_{i+2}, X'_{i+3}, \dots, X'_{m})\\
		&(Y_1, Y_{2}, \dots, Y_{i-1}, Y_i, C'_i \oplus C''_i, Y_{i+2}, Y_{i+3}, \dots, Y_m)
		\end{align*}
		also fulfil the conditions (S1) - (S3).
		
		Note that only the $i$th and $(i+1)$th elements of these sequences differ from the elements in the original sequences. 
		
		We can easily generalize this to the following: If $i < j_1 < j_2 < \dots < j_l$, there is an irreducible morphism $X_{j_k} \rightarrowtail X_{j_{k+1}}$ for all $1 \le k \le l$ and $X_{j_1} \mid X_i$, then there are two sequences with $X'_m$ and $Y_m$ as their $m$-th elements that together with
		\begin{align} \label{EqNewSeq}
		\begin{split}
		(X_1, \dots, X_{i-1}, X_{j_1}, X_{j_2}, \dots, X_{j_l}, X_i, X_{i+1}, \dots X_{j_1 - 1}, X_{j_1+1},\dots&\\
		\dots, X_{j_2 - 1}, X_{j_2+1}, \dots, X_{j_l-1}, X_{j_l+1}, \dots, X_l&)
		\end{split}
		\end{align}
		fulfil (S1) - (S3).
		
		Furthermore, there are sequences of modules that fulfil the conditions (S1) - (S3) with $X_m, X'_m, Y_m$ as their $m$-th elements so that $X''_i = X'_i$ for all $1 \le i \le m$: 
		
		By Definition \ref{DefBijection}, if there is a morphism $X_i \rightarrow X_j$, then $X_j < X_i$. So we can use the above to get sequences that fulfils (S1) -(S3) with $X_m, X'_m, Y_m$ as their $m$-th elements so that $X_1 \ge X_2 \ge \dots \ge X_{m-1}$. Then $X > \tau^{-1} I_i$ for all $X \mid X'_1$ and $X \mid Z_j$ with $j \le i$ and thus	$X''_i = X'_i$ for all $1 \le i < m$.
	\end{r1}
	
	\begin{proof}[Proof of Lemma \ref{LemExSeq}]
		We prove the lemma inductively. By Remark \ref{RemObs}, it is sufficient to prove the assertion for all sequences so that $X''_i = X'_i$ for all $1 \le i < m$.
		
		For these sequences, we additionally show the following: If there is an indecomposable direct summand $X$ of $X_m \oplus X'_m$ and $\tau X_i$ of $Y_m$ so that an irreducible morphism $X \rightarrow \tau X_i$ exists, then one of the following holds: 
		\begin{enumerate}[(a)]
			\item There is a direct summand $X' \cong X$ of $X_m \oplus X'_m$ so that the component $X' \rightarrow Y_m$ of $g_m$ is irreducible and $g_m(X') \subseteq \tau^{-1} X_i$
			
			\item Either $X \cong X_j$ for some $i < j < m$ or $X$ is isomorphic to a direct summand of $Y'_i$.
		\end{enumerate}
		
		If $m = 1$, the assertion is obvious by definition of the Auslander-Reiten sequence. 
		
		Now suppose that it holds for all series of modules of length $m \in \mathbb{N}$ or smaller. 
		We want to show that it also holds for sequences of length $m+1$ by applying Lemma \ref{LemARgen}. 
		
		To do this, we need to prove that there is an exact sequence of the form
		\begin{equation} \label{SeqAlg4}
		\xymatrixcolsep{3pc}\xymatrix{0 \ar[r] & M \ar[r] & X_m \oplus X'_m \ar^{{\left[\begin{smallmatrix}
					g_{m1} & g_{m2}\\ 0 & g_{m3}
					\end{smallmatrix}\right]}}[r] & Y'_m \oplus Y''_m \ar[r] & 0.}
		\end{equation}
		so that $g_{m1}$ is irreducible.
		
		Suppose that $Y'_m$ has some direct summands $Y'_{m1}, Y'_{m2}, \dots, Y'_{mk}$ and $g_m$ has a component 
		\begin{equation*}
		diag(g_{11}, g_{22}, \dots, g_{kk}):	X_m^k\rightarrow \bigoplus_{l = 1}^{k} Y'_{il}
		\end{equation*}
		where $g_{11}, g_{22}, \dots, g_{kk}$ are irreducible and $diag(g_{11}, g_{22}, \dots, g_{kk})$ is the diagonal matrix with entries $g_{11}, g_{22}, \dots, g_{kk}$. Then there is a copy of $X_m$ on which this restricts to
		\begin{equation*}
		\begin{bmatrix}
		g_{11}  \\ g_{22}\\ \vdots \\ g_{kk}
		\end{bmatrix}:	X_m \rightarrow \bigoplus_{l = 1}^{k} Y'_l,
		\end{equation*} 
		an irreducible morphism.
		
		By condition (S3) and since $Y''_1 = 0$, every indecomposable direct summand of $Y_m$ has the form $\tau^{-1} X_i$ for some $1 \le i < m$.
		
		If for all  $\tau^{-1} X_i \mid Y'_m$, there is some copy $X$ of $X_m$ so that the component $X \rightarrow \tau^{-1} Y_i$ of $g_m$ is irreducible and $g_m(X) \subseteq \tau^{-1} X_i$, then the above and the induction hypothesis mean that we can apply Lemma \ref{LemARgen}.
		
		Suppose that there is some  $\tau^{-1} X_i \mid Y'_m$, so that the above is not the case.
		
		Since $Y'_m \mid Z_m$, there is an irreducible morphism between $X_m$ and $\tau^{{-1}} X_i$. By the inductive hypothesis, one of the following holds:
		\begin{enumerate}[(a)]
			\item $X_m \cong X_j$ for some $i < j \le m$
			
			\item $X_m \mid Y'_i$.
		\end{enumerate}
		
		We show that there are sequences 
		\begin{align} \label{EqSeq(1)}
		\begin{split}
		&(X^{(1)}_1, X^{(1)}_2, \dots, X^{(1)}_m)\\
		&(X'^{(1)}_1, X'^{(1)}_2,\dots, X'^{(1)}_m)\\
		&(Y^{(1)}_1, Y^{(1)}_2, \dots, Y^{(1)}_m)
		\end{split}
		\end{align}
		that fulfil (S1)-(S3) and have $X_m$, $X'_m$, $Y_m$ as their $m$-th elements so that $\tau^{{-1}} X_i = \tau^{{-1}} X^{(1)}_{i'},$ but  $X_m \ncong X^{(1)}_j$ for all $i' < j \le m$ and $X_m \nmid Y'^{(1)}_{i'}$.
		
		Furthermore, we want to show that $X''^{(1)}_l = X'^{(1)}_l$ for all $1 \le l < m$. Since $X_m = X^{(1)}_m$, $X'_m = X'^{(1)}_m$, $Y_m = Y^{(1)}_m$, this is already clear for $i = m$.
		
		Obviously, we have $X_m \mid X_m \oplus  X'_m$, so either $X_m \mid X'_1$ or $X_m \mid Z'_k$ for some $1 \le k \le m$, $k \neq i$. 
		
		In the first case, (b) is not possible, since $Y_i$ and $X_i \oplus X'_i$ do not share direct summands. In case (a), Remark \ref{RemObs} yields a sequence $(X^{(2)}_1, \dots, X^{(2)}_m)$, where $X_j$ comes before $X_i$. 
		
		In the second case, we can get a new sequence where $X_k$ comes before $X_i$, since $Z'_i \mid X_m \oplus X'_m$ (otherwise, $\tau^{-1}X_i$ would not be a direct summand of $Y_m$). In case (b), this sequence is already the one we need; in case (a), we can again get a another sequence by Remark \ref{RemObs} where $X_j$ comes before $X_i$.
		
		If we call this new exact sequence $(X^{(2)}_1, \dots, X^{(2)}_m)$, then it is clear  by (\ref{EqNewSeq}) that $X''^{(2)}_l = X'^{(2)}_l$ holds for all $1 \le l < m$.
		
		Since there are only finitely many $j$ with $i' < j \le m$, we get sequences of the form (\ref{EqSeq(1)}) after finitely many steps.
		
		The inductive assumption gives us an exact sequence
		\begin{equation} \label{SeqAlg7}
		\xymatrix{0 \ar[r] & M \ar^-{f'_m}[r]& X_m \oplus X'_m \ar^-{g'_m}[r] & Y_m \ar[r] & 0}
		\end{equation}
		where the component $X_m \rightarrow Y_m$ of $g'_m$ is irreducible and $g'_m(X_m) \subseteq \tau^{-1} X_i$.
		
		If $\tau^{-1} X_i = Y_m$, then it is sufficient to look at the sequence (\ref{SeqAlg7}) instead of  
		\begin{equation} \label{SeqAlg5}
		\xymatrix{0 \ar[r] & M \ar^-{f_m}[r] & X_m \oplus X'_m \ar^-{g_m}[r] & Y_m \ar[r]& 0}.
		\end{equation} 
		
		If there is some $\tau^{-1} X_k$  so that $\tau^{-1} X_i \oplus \tau^{-1} X_k \mid Y_m$, then we can assume that $g_m$ induces an indecomposable morphism $X_m \rightarrow \tau^{-1} X_k$ and $g_m(X_m) \subset \tau^{-1} X_k$.
		
		So together (\ref{SeqAlg5}) and (\ref{SeqAlg7}) give a new exact sequence 
		\begin{equation*}
		\xymatrix{0 \ar[r] & M \ar^-{f''_m}[r]& X_m \oplus X'_m \ar^-{g''_m}[r] & Y_m \ar[r] & 0},
		\end{equation*}
		where the induced morphisms $X_m \rightarrow \tau^{-1} X_i$ and $X_m \rightarrow \tau X_k$ of $g''_m$ are irreducible and $g_m(X_m) \subset \tau^{-1} X_i \oplus \tau^{-1} X_k$. 
		
		Inductively, there is an exact sequence of the form (\ref{SeqAlg4}), where $g_{m1}$ is irreducible and we can use Lemma \ref{LemARgen} to get an exact sequence
		\begin{equation*}
		\xymatrix{0 \ar[r] & M \ar^-{f_{m+1}}[r] & X_{m+1} \oplus X'_{m+1} \ar^-{g_{m+1}}[r] & Y_{m+1} \ar[r]& 0}.
		\end{equation*}
		If there is a monomorphism $M \rightarrowtail U$, then it factors through $f_{m+1}$.
		
		This gives us not only the assertion of the lemma but also the additional assumptions we have made:
		
		Let 
		\begin{equation*}
		g_m = \begin{bmatrix}
		g'_{m1} & g'_{m2}
		\end{bmatrix}: Y'_m \oplus Z'_m \rightarrow \tau^{-1} X_m
		\end{equation*}
		be the epimorphism of the AR-sequence. By (\ref{SeqAlg4}) and Lemma \ref{LemARgen}, we get 
		\begin{equation*}
		g_{m+1} = \begin{bmatrix}
		g'_{m2} & g'_{m1} g_{m2} \\
		0 & g_{m3}
		\end{bmatrix}: Z'_{m} \oplus X'_m \rightarrow \tau^{-1} X_m \oplus Y''_m.
		\end{equation*}
		Let $X$ be a direct summand of $Z'_m \oplus X'_m$ and $\tau^{-1} X_i$ a direct summand of $\tau^{-1} X_m \oplus Y''_m$ so that there is an irreducible morphism $X \rightarrow \tau^{-1} X_i$.
		
		If $i = m$, then $X$ is a direct summand of $Z_m$. If it is also a direct summand of $Z'_m$, then $g_{m+1}(X) = \left[\begin{smallmatrix}
		g'_{m2}\\
		0
		\end{smallmatrix}\right](X)$ and (a) holds. Otherwise, $X$ is a direct summand of $Y'_m$.
		
		If $i \neq m$ and $X$ is a direct summand of $X'_m$, then either (b) holds or
		\begin{equation*}
		\begin{bmatrix}
		g_{m2} \\g_{m3}
		\end{bmatrix}(X'_m) \subset \tau^{-1} X_i.
		\end{equation*}
		Thus $g_{m2}(X'_m) = 0$ and $g_{m+1}(X) = \left[\begin{smallmatrix}
		0\\
		g_{m3}
		\end{smallmatrix}\right](X)$. So (a) holds.
		
		Finally, suppose that $i \neq m$ and $X$ is not a direct summand of $X'_m$. Because  of the irreducible morphism between $X$ and $\tau^{-1} X_i$, the former is a direct summand of $Z_i$. By (S3), either it is a direct summand of $Y'_i$ or of $X_j$ for some $i < j < m$.
	\end{proof}
	
	A perhaps simpler way to interpret the sequences of modules used in the lemma above is the following:
	
	\begin{r1}
		Suppose that $X_i, X'_i, Y_i$ are the $i$-th elements of sequences that fulfil (S1) - (S3). Then $X_{i+1} \oplus X'_{i+1}$ are defined by taking the exact sequence
		\begin{equation*}
		\xymatrix{0 \ar[r] & M \ar[r] & X_i \oplus X'_i \ar[r] & Y_i \ar[r]& 0}.
		\end{equation*}
		and the Auslander-Reiten sequence
		\begin{equation*}
		\xymatrix{0 \ar[r] & X_i \ar[r] & Z_i \ar[r] & \tau^{-1} X_i \ar[r]& 0}.
		\end{equation*}
		We can add these sequences together and get
		\begin{equation*}
		\xymatrix{0 \ar[r] &M \oplus  X_i \ar[r] & X_i \oplus X'_i \oplus Z_i \ar[r] & Y_i \oplus \tau^{-1} X_i \ar[r]& 0}.
		\end{equation*}
		Then $X_i$ is the maximal module that is a direct summand of both the first and the second term: we still get an exact sequence if we delete it in both terms:
		\begin{equation*}
		\xymatrix{0 \ar[r] &M \ar[r] &  X'_i \oplus Z_i \ar[r] & Y_i \oplus \tau^{-1} X_i \ar[r]& 0}.
		\end{equation*}
		The same holds for $G_i$, the maximal module that is both a direct summand of the middle term and the last term. Deleting this in both terms gives us an exact sequence
		\begin{equation*}
		\xymatrix{0 \ar[r] & M \ar[r] & X_{i+1} \oplus X'_{i+1} \ar[r] & Y_{i+1} \ar[r]& 0}.
		\end{equation*}
	\end{r1}

	These modules have some interesting properties:
	
	\begin{c1} \label{CorMono}
		If there is a monomorphism $h: M \rightarrowtail U$, then for sequences of modules 
		\begin{align*}
		&(X_1, X_2, \dots, X_m)\\
		&(X'_1, X'_2,\dots, X'_m)\\
		&(Y_1, Y_2, \dots, Y_m),
		\end{align*}
		which fulfil (S1) - (S3), there is a monomorphism $X_i \oplus X'_i \rightarrowtail Y_i \oplus U$ for all $1 \le i < m$.
		
		Thus, every injective direct summand of $X_i \oplus X'_i$ is a direct summand of $U$.
	\end{c1}
	
	\begin{proof}
		By Lemma \ref{LemExSeq}, there is an exact sequence
		\begin{equation*}
		\xymatrix{0 \ar[r] & M \ar^-{f_i}[r] & X_i \oplus X'_i \ar^-{g_i}[r] & Y_i \ar[r]& 0}
		\end{equation*}
		for $1 \le i \le m$ so that $h$ factors through $f_i$. Thus, there is some morphism $h_i$ with $h = h_if_i$. So $h_i$ is a monomorphism on $\Im f_i$. Since $\Ker g_i = \Im f_i$, the morphism
		\begin{equation*}
		\begin{bmatrix}
		g_i \\ h_i
		\end{bmatrix}: X_i \oplus X'_i \rightarrowtail Y_i \oplus U
		\end{equation*}
		is a monomorphism.
		
		So every injective direct summand $I$ of $X_i \oplus X'_i$ is a direct summand of $Y_i \oplus U$. Since $X_i \oplus X'_i$ and $Y_i$ do not share any direct summands, $I$ is even a direct summand of $U$.
	\end{proof}
	
	We can use the following lemma to show that there is an algorithm that, for given indecomposable, preinjective module $M$ constructs all $U$ with $M \rightarrowtail U$.
	
	\begin{l1} \label{LemARNumbers}
		Suppose there is an irreducible morphism between $(s, i) = \tau^{s} I_i$ and $(t, j) = \tau^{t} I_j$.
		Then either $s=t$ and $i > j$ or $s = t+1$ and $i < j$.
	\end{l1}
	
	\begin{proof}
		By \cite{ARS}, Chapter VIII, Proposition 1.2 and Lemma 1.8, 
		\begin{equation*}
		s-1 \le  t \le s.
		\end{equation*}
		Furthermore, we have ordered the injective modules so that $\Hom(I_i, I_j) = 0$ if $i < j$. 
		
		By \cite{ARS} Chapter VIII, Corollary 4.2, we get 
		$\Hom(\tau^s I_i, \tau^s I_j) = 0$ if $i < j$ and by  \cite{ARS} Chapter VIII, Corollary 4.3,, there is an irreducible morphism $\tau^s I_i \rightarrow \tau^{s-1} I_j$ if and only if there is an irreducible morphism $\tau^{s} I_j \rightarrow \tau^s I_i$.
		
		So if $s =t$, then $i > j$ and if $s = t+1$ then $i < j$. 
	\end{proof}
	
	\begin{p1} \label{Alg}
		Let $A$ be a hereditary Artin algebra with $M \in \mod A$ indecomposable and preinjective. Let $U \in \mod A$, so that $M$ is not a direct summand of $U$. There is a monomorphism $M \rightarrowtail U$ if and only if for some $m \in \mathbb{N}$ there are three sequences of modules
		\begin{align*}
		&(X_1, X_2, \dots, X_m)\\
		&(X'_1, X'_2,\dots, X'_m)\\
		&(Y_1, Y_2, \dots, Y_m)
		\end{align*}
		that fulfil the conditions (S1) - (S3) and furthermore 
		\begin{enumerate}
			\item[(S4)] If for some $1 \le i \le m$ the module $X_i \oplus X'_i$ has an injective direct summand $I$, then $I \mid U$.
			
			\item[(S5)] $X_m \oplus X'_m$ is a direct summand of $U$.
		\end{enumerate}
	\end{p1}
	
	\begin{proof} To prove this, we use Lemma \ref{LemExSeq}: since the sequences fulfil (S1)-(S3), there are exact sequences of the form
		\begin{equation*}
		\xymatrix{0 \ar[r] & M \ar^-{f_i}[r] & X_i \oplus X'_i \ar^-{g_i}[r] & Y_i \ar[r]& 0}
		\end{equation*}
		for all $1 \le i \le m$. If a monomorphism $M \rightarrowtail U$ exists, it factors through $f_i$ for all $1 \le i \le m$. 
		
		Thus one direction is obvious: if such sequences of modules exist, $f_m: M \rightarrowtail U$ is a monomorphism. 
		
		On the other hand, suppose that no series of modules fulfil (S1) - (S5). 
		
		If $M$ is injective, then it cannot be a submodule of $U$. Otherwise, there are series of modules that fulfil (S1) - (S3), since there is an AR-sequence that starts in $M$ and we can set $m = 1$. 
		
		If (S4) is not fulfilled, then $M$ cannot be a submodule of $U$ by Corollary \ref{CorMono}. 
		
		Otherwise, there is some non-injective $X_{m+1}$ and some $\alpha_{m+1} \in \mathbb{N}$ so that $X_{m+1}^{\alpha_{i+1}} \mid X_{m} \oplus X'_{m}$, but $X_{m+1}^{\alpha_{i+1}} \nmid U$. So we can extend the sequences of modules to 
		\begin{align*}
		&(X_1, X_2, \dots X_m, X_{m+1})\\
		&(X'_1, X'_2, \dots X'_m, X'_{m+1})\\
		&(Y_1, Y_2, \dots Y_m, Y_{m+1})
		\end{align*}
		so that these series fulfil (S1) - (S3). If these sequences fulfil (S4), we can extend them again to sequences of length $m +2$. 
		
		We have $M = (r, i)$ for some $r \in \mathbb{N}$ and $1 \le i \le n$. Every indecomposable direct summand of $X_1 \oplus X'_1$ is of the form $(r', j) < (r, i)$ for some $r' \in \mathbb{N}_0$ and $1 \le j \le n$. Furthermore, if $X_1 = (r', j)$, then every direct summand of $Z'_1$ is of the form $(r'', k) < (r', j)$, and analogously for $X_2, X_3, \dots$. 
		
		So after finitely many steps, either we find sequences that do not fulfil (S4), or there is some $m'$ so that every direct summand of $X_m \oplus X'_m$ is injective. If (S4) is still fulfilled, then (S5) is also fulfilled, a contradiction to our assumption.
	\end{proof}
	
	The proof of Proposition \ref{Alg} shows the following:
	
	\begin{c1} \label{CorAlg}
		Let $M$ and $U$ be preinjective modules over $A$. If $M \subset U$, then all sequences of modules that fulfil (S1) - (S3) can be extended to sequences of modules that fulfil (S1) - (S5).
		
		If $M \not\subset U$, then all sequences of modules that fulfil (S1) - (S3) can be extended to sequences that fulfil (S1) - (S3) so that $X_m \oplus X'_m$ has an injective direct summand that is not a direct summand of $U$.
	\end{c1}
	
	\begin{r1} \label{RemAlg}
		By Corollary \ref{CorAlg}, we can use the proposition as an algorithm that finds out for given indecomposable, preinjective $M$ and modules $U$, if there is a monomorphism $M \rightarrowtail U$. Alternatively, we can use it to construct all $U$ with $M \subset U$.
	\end{r1}
	
	Note that it is very simple to generalize this for arbitrary preinjective $M$:
	
	\begin{c1} \label{CorRemAlg}
		Let $M$ be a preinjective module so that $M = \bigoplus_{i=1}^m M_i$ with $M_i$ indecomposable. Let $U$ be some module in $\mod A$. Denote the middle term of the Auslander-Reiten sequence that starts in $M_i$ by $N_i$. 
		
		Furthermore, order $M_1, \dots, M_m$ so that there is some $0 \le k \le m$ with $M_i \mid U$ if and only if $i \le k$.
		
		Suppose that the sequences of modules
		\begin{align*}
		&(X_1, X_2, \dots, X_m)\\
		&(X'_1, X'_2,\dots, X'_m)\\
		&(Y_1, Y_2, \dots, Y_m)
		\end{align*}
		fulfil (S2), (S3) and 
		\begin{enumerate}[(S'1)]
			\item We have 
			\begin{equation*}
			X_1 \oplus X'_1 = \bigoplus_{i = 1}^k M_i \oplus \bigoplus_{i = k+1}^m N_i
			\end{equation*}
			and
			\begin{equation*}
			Y_1 = \bigoplus_{i = k+1}^m \tau^{-1} M_i.
			\end{equation*}
		\end{enumerate}
		Then for all $1 \le i \le m$, there is an exact sequence
		\begin{equation*}
		\xymatrix{0 \ar[r] & M \ar^-{f_i}[r] & X_i \oplus X'_i \ar^-{g_i}[r] & Y_i \ar[r]& 0}.
		\end{equation*}
		There is a monomorphism $M \rightarrowtail U$ if and only if there is some $m' > m$ and modules $X_{m+1}, \dots X_{m'}$, $X'_{m+1}, \dots X'_{m'}$, $Y_{m+1}, \dots Y_{m'}$ so that the sequences	
		\begin{align*}
		&(X_1, X_2, \dots, X_{m'})\\
		&(X'_1, X'_2,\dots, X'_{m'})\\
		&(Y_1, Y_2, \dots, Y_{m'})
		\end{align*}
		fulfil (S'1) and (S2) - (S5).
		
		Furthermore, if a monomorphism $M \rightarrowtail U$ exists, then it factors through all $f_i$.
	\end{c1}

	\begin{e1}
		Take $A$ as in Example \ref{ExWeyl}. A part of the preinjective component of the AR-quiver of $A$ is:
		\[
		\xymatrix {
			\dots \ar[r] \ar[dr]& \tau I_3 \ar[r] \ar[dr]&\tau I_1 \ar[r] \ar[dr]&I_3 \ar[r] \ar[dr]&I_1\\
			\dots \ar[r] \ar[ur] & \tau I_4 \ar[r] \ar[ur]&\tau I_2\ar[r] \ar[ur]&I_4 \ar[r] \ar[ur] &I_2
		}
		\]
		Suppose that we want to know whether $M = \tau I_3$ is a submodule of, say, $U = I_2 \oplus I_3\oplus I_4$. 
		
		Then by (S1), $X_1 \oplus X'_1 = \tau I_1 \oplus \tau I_2$ and $Y_1 = I_3$. Since neither $\tau I_1$ nor $\tau I_2$ is a direct summand of $U$, we arbitrarily set $X_1 := \tau I_1$.
		
		The AR-sequence
		\[
		\xymatrix{0 \ar[r] & \tau I_1 \ar[r] & I_3 \oplus I_4 \ar[r] & I_1 \ar[r] & 0}
		\]
		and (S3) show that $X_2 \oplus X'_2 = \tau I_2 \oplus I_4$ and $Y_2 = I_1$. Since $I_4$ is a direct summand of $U$, we set $X_2 := \tau I_2$ to fulfil (S2).
		Using the AR-sequence
		\[
		\xymatrix{0 \ar[r] & \tau I_2 \ar[r] & I_3 \oplus I_4 \ar[r] & I_2 \ar[r] & 0}
		\]
		we get $X_3 \oplus X'_3 = I_3 \oplus I^2_4$ and $Y_3 = I_1 \oplus I_2$. Since $I^2_4$ is injective, but not a direct summand of $U$, the condition (S4) is not fulfilled and there is no monomorphism between $M$ and $U$. 
	\end{e1}
	
	We have one more lemma:

	\begin{l1} \label{LemAlgInd}
		Let $M$ be an indecomposable, preinjective module and $U \in \mod A$. If the sequences  
		\begin{align} \label{Sequences}
		\begin{split}
		&(X_1, X_2, \dots, X_m)\\
		&(X'_1, X'_2,\dots, X'_m)\\
		&(Y_1, Y_2, \dots, Y_m)
		\end{split}
		\end{align}
		fulfil (S1) - (S3), then for every $1 \le i \le m$, there is an exact sequence
		\begin{equation} \label{Sequence}
		\xymatrix{0 \ar[r] & X_i \oplus X'_i \ar[r] & Y_i \oplus X_m \oplus X'_m \ar[r] & Y_m \ar[r] & 0}.
		\end{equation}
		Furthermore, if there is an exact sequence 
		\begin{equation} \label{Sequence2}
		\xymatrix{0 \ar[r] & X_i \oplus X'_i \ar[r] & Y_i \oplus U \ar[r] & Z \ar[r] & 0}
		\end{equation}
		then there is also an exact sequence
		\begin{equation} \label{Sequence3}
		\xymatrix{0 \ar[r] & M_0 \ar[r] & U \ar[r] & Z \ar[r] & 0}.
		\end{equation}
	\end{l1}
	
	\begin{proof}
		Let $Y$ be the maximal module so that $Y \mid Y_i$ and $Y \mid Y_m$. Furthermore, suppose that $Y_i = Y \oplus Y'$ and $Y_m = Y \oplus Y''$.
		
		We use Corollary \ref{CorRemAlg} on $X_i \oplus X'_i$ and $Y \oplus U$. Take $i < j_1 < j_2 \dots < j_l$ so that $X_{j_k}$ are those modules in the sequence $(X_{i+1}, \dots X_{m})$ which are already a direct summand of $X_i \oplus X'_i$. Then 
		\begin{equation*}
		(X_{m+1}, \dots, X_{j_1-1}, X_{j_1+1}, \dots, X_{j_2-1}, X_{j_2+1}, \dots, X_{j_l-1}, X_{j_l+1}, \dots, X_m)
		\end{equation*}
		is part of a triple of sequences that fulfil (S'1) and (S2) - (S5) with respect to $X_i \oplus X'_i$ and $Y' \oplus U$.
		
		So the same construction that yields an exact sequence
		\begin{equation*}
		\xymatrix{0 \ar[r] & M \ar[r] & X_m \oplus X'_m \ar[r] & Y_m \ar[r] & 0},
		\end{equation*}
		also gives an exact sequence
		\begin{equation*}
		\xymatrix{0 \ar[r] & X_i \oplus X'_i \ar[r] & Y' \oplus X_m \oplus X'_m \ar[r] & Y'' \ar[r] & 0},
		\end{equation*}
		when used on $X_i \oplus X'_i$ and $Y' \oplus U$ instead of $M$ and $U$. Adding $Y$ to both the middle and the last term gives (\ref{Sequence}).
		
		The exact sequence (\ref{Sequence2}) is given by a sequence of modules that fulfil (S'1) and (S2) - (S5). Together with the sequences (\ref{Sequences}), this yields the exact sequence \ref{Sequence3}.
	\end{proof}

	\section{Preinjective modules and the Weyl group} \label{SecConnection}
	
	In this section we connect our results about preinjective modules with the relations of the Weyl group.
	
	First we give a Lemma that shows the connection between the AR-sequences and the relations:
	
	\begin{l1} \label{LemARWeyl}
		Fix two integers $1 \le i, j \le n$. Set
		\begin{align*}
		\alpha &:= \max\{v ~|~ \exists (s,i), (t,j): \text{ there is an irreducible morphism } (s,i) \rightarrow (t, j)^{v} \}\\
		\beta &:= \max\{v ~|~ \exists (s,i), (t,j): \text{ there is an irreducible morphism } (t,j) \rightarrow (s, i)^{v} \}
		\end{align*}
		Let $(s_i s_j)^{m_{ij}}$be the defining relation of the Weyl group as in Lemma \ref{LemCoxeter}.
		
		Then the value of $m_{ij}$ depends on $\alpha \beta$ in the following way:
		\begin{equation*}
		\begin{array}{l|llllll} \alpha \beta& 0 & 1 & 2 & 3 & \ge 4\\ \hline
		m_{ij} & 2& 3 & 4 & 6 & \infty.
		\end{array}
		\end{equation*}
	\end{l1}
	
	\begin{proof}
		By \cite{ARS} VIII, Corollaries 4.2 and 4.3, the integers $\alpha$ and $\beta$ do not depend on $s$. Let $(c_{ij})_{nn}$ be the Cartan matrix. Then by \cite{ARS}, p. 267-268,, either $\alpha = c_{ij}$ and $\beta = c_{ji}$ or $\beta = c_{ij}$ and $\alpha = c_{ji}$. Lemma \ref{LemCoxeter} gives the stated values for $m_{ij}$.
	\end{proof}
	
	Now we define a recursion that plays a fundamental role in the proof of the bijection:
	
	\begin{d1} \label{DefRec}
		For given $\alpha, \beta \in \mathbb{N}$ define a recursion formula by
		\begin{align*}
		&E(0) = 1\\
		&E(1) = \alpha\\
		&E(2m) = \beta E(2m-1) - E(2m - 2)\\
		&E(2m+1) = \alpha E(2m) - E(2m-1)		
		\end{align*}
		for all $m \in \mathbb{N}$.
	\end{d1}
	
	This recursion is directly linked to the Weyl group:

	\begin{l1} \label{LemRec}
		Let $\alpha$, $\beta$ be as in Lemma \ref{LemARWeyl}. 
		Then 
		\begin{equation*}
		E(m) = 0 \Longleftrightarrow m \ge m_{ij} -1.
		\end{equation*}
	\end{l1}
	
	\begin{proof}
		If $\alpha \beta < 4$, then we get the following values for $m \le 6$:
		
		\begin{equation} \label{Table}
		\begin{array}{l|l|lllllll} \alpha \beta & m_{ij}& E(0) & E(1) & E(2) & E(3) & E(4) & E(5) & E(6)\\ \hline
		0 & 2& 1 & 0 & 0 & 0 & 0 & 0 & 0 \\
		1 & 3&1 & 1 & 0 & 0 & 0 & 0 & 0\\
		2 & 4&1 & \alpha & 1 & 0 & 0 & 0 & 0\\
		3 &6& 1 & \alpha & 2 & \alpha & 1 & 0 & 0\\ 
		\end{array}
		\end{equation}
		
		Obviously, $E(m) = 0$ if $\alpha \beta < 4$ and $m > 6$.
		
		If $\alpha \beta \ge 4$, then $m_{ij} = \infty$ by Lemma \ref{LemARWeyl} and we need to show that $E(m) \neq 0$ for all $m \in \mathbb{N}$.
		
		Since $E(2) = \alpha \beta - 1 > E(0) > 0$, we get inductively for $m > 1$:
		\begin{align*}
		E(2m) &= \beta E(2m-1) - E(2m-2)\\
		&= (\alpha \beta -1)E(2m-2)- \beta E(2m-3)\\
		&= (\alpha \beta -1)E(2m-2)- E(2m-2) -E(2m-4)\\
		&> (\alpha \beta -3)E(2m-2)\\
		&\ge E(2m-2).
		\end{align*}
		The proof that $E(2m +1) > E(2m-1)>0$ is completely analogous.
	\end{proof}
	
	Next, we need some notation:
	
	\begin{d1} \label{DefMs}
		Fix $s \in \mathbb{N}_0$ and $1 \le i \neq j \le n$ and let $M_0 := \tau^s I_i$. If $s \ge 1$ or $j < i$, let $t$ be the integer with $(s-1, i) < (t, j) < (s, i)$. Denote $M_1 := \tau^t I_j$, $M_2 := \tau^{s-1} I_i$, $M_3 := \tau^{t-1} I_j$, $\dots$.
	\end{d1}
	
	The following lemma is a key part in the proof that there is a bijection between cofinite, submodule closed subcategories and the elements of the Weyl group:
	
	\begin{l1} \label{LemRecSeq}
		Let $U$ be a module so that $M_k \nmid U$ for all $M_k \neq 0$ with $0 \le k \le m_{ij}-1$.  Then for all $m \ge 1$ with $M_{m+1}\neq 0$ and $E(m-1) > 0$, there are series of modules that fulfil (S1) - (S3) and yield exact sequences
		\begin{equation} \label{Seq}
		\xymatrix{0 \ar[r] & M_0 \ar^-{f_m}[r] & M_m^{E(m)} \oplus U_m \ar[r] & M_{m+1}^{E(m-1)} \ar[r] & 0}
		\end{equation}
		so that no $M_k \nmid U_m$ for all $M_k \neq 0$ with $0 \le k \le m_{ij}-1$.
		
		If there is a monomorphism $M_0 \rightarrowtail U$, then it factors through $f_m$ for all $m$.
	\end{l1}
	
	\begin{proof}
		If $M_0$ is injective, there is nothing to show. So we can assume that an AR-sequence starts in $M_0$.
		
		Let $\alpha, \beta$ be as in Lemma \ref{LemARWeyl}. Then there are modules $M', N'$ so that
		\begin{equation} \label{SeqM}
		\xymatrix{0 \ar[r] & M_0 \ar[r] & M_1^\alpha \oplus N' \ar[r] & M_2 \ar[r] & 0}
		\end{equation}
		and
		\begin{equation*}
		\xymatrix{0 \ar[r] & M_1 \ar[r] & M_2^\beta \oplus M'  \ar[r] & M_3 \ar[r] & 0}
		\end{equation*}
		are AR-sequences. Note that by \cite{ARS} VIII, Corollaries 4.2 and 4.3, for all non-injective $M_{2m-1}$, $m \in \mathbb{N}_0$, there are AR-sequences of the form
		\begin{equation} \label{ARtau1}
		\xymatrix{0 \ar[r] & M_{2m-1} \ar[r] & M_{2m}^\beta \oplus \tau^m M'  \ar[r] & M_{2m+1} \ar[r] & 0}.
		\end{equation}
		For all non-injective $M_{2m}$, $m \in \mathbb{N}$ they are of the form
		
		\begin{equation*}
		\xymatrix{0 \ar[r] & M_{2m} \ar[r] & M_{2m+1}^\alpha \oplus \tau^m N' \ar[r] & M_{2m+2} \ar[r] & 0}.
		\end{equation*}
		
		If we set $U_1 := N'$, the AR-sequence that starts in $M_0$ is the exact sequence
		\begin{equation*}
		\xymatrix{0 \ar[r] & M_0 \ar[r] & M_1^{E(1)} \oplus U_1 \ar[r] & M_{2}^{E(0)} \ar[r] & 0}.
		\end{equation*}
		If $M_1$ is injective, then the proof is complete. So we can assume that an AR-sequence starts in $M_1$ and use Lemma \ref{LemExSeq}. 
		
		Since $M_1 \nmid U$, we set
		\begin{equation*}
		X_1 := X_2 := X_3 := \dots := X_{E(1)} := M_1, 
		\end{equation*}
		Then we get sequences that fulfil conditions (S1) - (S3) by setting
		\begin{align*}
		X'_1 &:= M^{E(1)-1}_1 \oplus N',\\
		X'_2 &:= M^{E(1)-2}_1 \oplus M^{\beta-1}_2 \oplus N' \oplus M',\\
		&\dots\\
		X'_{E(1)} &:= M_2^{(E(1)-1)\beta-1} \oplus N' \oplus (M')^{E(1)-1}\\
		\intertext{and}
		Y_1 &:= M_2,\\ Y_2 &:= M_3,\\ Y_3 &:= M^2_3,\\ \dots\\ Y_{E(1)} &:= M^{E(1)-1}_3
		\end{align*}
		
		We get 
		\begin{equation*}
		X_{E(1)+1} \oplus X'_{E(1)+1} = M_{2}^{E(1)\beta - 1} \oplus N' \oplus (M')^{E(1)}
		\end{equation*} 
		and $Y_{E(1)+1} = M^{E(1)}_3$. Thus by Lemma \ref{LemExSeq} there is some $f_1$ so that the following sequence is exact:
		\begin{equation*}
		\xymatrix{0 \ar[r] & M_0 \ar^-{f_1}[r] & M_{2}^{E(1)\beta - 1} \oplus N' \oplus (M')^{E(1)} \ar[r] & M^{E(1)}_3 \ar[r] & 0}.
		\end{equation*}
		If there is a monomorphism $M \rightarrowtail U$, then it factors through $f_1$. 
		
		Since $U_2 := N' \oplus (M')^{E(1)}$, we can write the exact sequence as
		\begin{equation*}
		\xymatrix{0 \ar[r] & M_0 \ar^-{f_1}[r] & M_{2}^{E(2)} \oplus U_2 \ar[r] & M^{E(1)}_3 \ar[r] & 0}.
		\end{equation*}
		We show the rest inductively: Suppose that
		\begin{equation*}
		\xymatrix{0 \ar[r] & M_0 \ar[r] & M_{2m-1}^{E(2m-1)} \oplus U_{2m-1} \ar[r] & M_{2m}^{E(2m-2)} \ar[r] & 0}
		\end{equation*}
		is an exact sequence and $E(2m-1) \neq 0$. Furthermore, suppose that this exact sequence is yielded by sequences of modules of the length $m'-1$. Then  
		\begin{equation*}
		X_{m'} := M_{2m-1}, X'_{m'+1} := M_{2m-1}^{E(2m-1)-1} \oplus U_{2m-1}, Y_{m'} := M_{2m}^{E(2m-2)}
		\end{equation*}
		are elements of sequences that fulfil the condition (S1) - (S3) of Lemma \ref{LemExSeq}.
		
		If $M_{2m-1}$ is injective, then $M_{2m+1} = \tau M_{2m-1} = 0$ and there is nothing to prove. If $M_{2m-1}$ is not injective, then the AR-sequence (\ref{ARtau1}) exists. As above, we set
		\begin{equation*}
		X_{m'+1} := \dots := X_{m'+ E(2m-1)-1} := M_{2m-1}.
		\end{equation*}
		This determines $X'_{m'+1}, \dots, X'_{m'+ E(2m-1)-1}$ and $Y_{m'+1}, \dots, Y_{m'+ E(2m-1)-1}$ completely.
		
		Since $E(2m) = \beta E(2m-1)- E(2m-2)$, we get 
		\begin{align*}
		X_{m'+ E(2m-1)} \oplus X'_{m'+ E(2m-1)} &= M_{2m}^{E(2m)} \oplus U_{2m-1} \oplus (\tau^{m-1} M')^{E(2m-1)}\\
		Y_{m'+ E(2m-1)} &= M_{2m+1}^{E(2m-1)}
		\end{align*}
		Together with $U_{2m} :=U_{2m-1} \oplus (\tau^{m-1} M')^{E(2m-1)}$, this yields an exact sequence 
		\begin{equation*}
		\xymatrix{0 \ar[r] & M_0 \ar^-{f_{2m}}[r] & M_{2m}^{E(2m)} \oplus U_{2m} \ar[r] & M_{2m+1}^{E(2m-1)} \ar[r] & 0}
		\end{equation*}
		for some $f_{2m}$. By Lemma \ref{LemExSeq}, $M \rightarrowtail U$ factors through $f_{2m}$.
		
		Analogously, we can construct 
		\begin{equation*}
		\xymatrix{0 \ar[r] & M_0 \ar^-{f_{2m+1}}[r] & M_{2m+1}^{E(2m+1)} \oplus U_{2m} \ar[r] & M_{2m}^{E(2m+2)} \ar[r] & 0}
		\end{equation*}
		if $E(2m) \neq 0$ and $M_{2m+2} \neq 0$.
	\end{proof}
	
	\begin{c1} \label{CorRecSeq}
		Let $U_m$ be as in Lemma \ref{LemRecSeq}. If $(r, l) \mid U_m$, then $M_0 > (r, l) > M_{m+1}$.
		
		If $M_0 > (r, l) > M_1$, then $(r, l) \mid U_m$ if and only if $m_{il} \ge 3$. 
		
		If $M_1 > (r, l) > M_m$, then $(r, l)\mid U_m$ if and only if $m_{il} + m_{jl} \ge 5$.
		
		If $M_m > (r, l) > M_{m+1}$ and $m$ is even, then $(r, l)\mid U_m$ if and only if $m_{il} \ge 3$. If $m$ is odd, then $(r, l)\mid U_m$ if and only if $m_{jl} \ge 3$.
	\end{c1}
	
	\begin{proof}
		This is obvious from the proof of Lemma \ref{LemRecSeq}.
	\end{proof}
	
	\begin{r1} \label{RemRecSec}
		Note that $m_{ij} = m_{ji}$. If we fix $s, t$ as in Definition \ref{DefMs}, we can set $M'_0 := \tau^t I_j$, $M'_1 := \tau^{s-1} I_i$, $M'_2 := \tau^{t-1} I_j$, $M'_3 := \tau^{s-2} I_i$, $\dots$ and
		\begin{align*}
		E'(0) &= 1\\
		E'(1) &= \beta\\
		E'(2m) &= \alpha E(2m-1) - E(2m - 2)\\
		E'(2m+1)&= \beta E(2m) - E(2m-1).		
		\end{align*}
		With this definition, we get analogous results to \ref{LemRec}, \ref{LemRecSeq} and \ref{CorRecSeq}.
	\end{r1}

	\section{Preliminaries for the main theorem} \label{SecProof}
	
	This section collects some preliminaries which are necessary to prove that there is a bijection between the Weyl group and the cofinite, submodule closed subcategories: First, we show that every cofinite submodule closed subcategory is of the form $\mathcal{C}_{\underline{w}}$ for some word $\underline{w}$.
	
	Then we will prove an auxiliary result that will make the inductions in the next section possible.
	
	\begin{l1} \label{LemWord}
		If a cofinite, full additive subcategory $\mathcal{C}$ of $\mod A$ is submodule closed, then there is a word $\underline{w}$ over $S =\{s_1, s_2, \dots s_n\}$ with $\mathcal{C} = \mathcal{C}_{\underline{w}}$.
	\end{l1}
	
	\begin{proof}
		By Lemma \ref{LemPreinj},  
		\begin{equation*}
		\ind A \setminus \mathcal{C} = \ind \mathcal{I} \setminus \mathcal{C} =: \{(r_1, i_1), (r_2, i_2), \dots (i_m, r_m)\}.
		\end{equation*}
		for some $m \in \mathbb{N}$ and modules $(r_1, i_1) < (r_2, i_2) < \dots < (r_m, i_m)$.
		
		Suppose that for all words $\underline{w}$ over $S$
		\begin{equation*}
		\rho(\underline{w})\neq (r_1, i_1)(r_2, i_2) \dots (r_m, i_m).
		\end{equation*} 
		By Definition \ref{DefWord}, either $r_1 > 0$ or there is some $1 \le j \le m-1$ so that
		\begin{equation*}
		(r_j, i_j) < (r_{j+1}, i_{j+1}-1).
		\end{equation*}
		In the first case,  $\mathcal{C}$ contains the middle term of the AR-sequence that starts in $(r_1, i_1)$ by Lemma \ref{LemARNumbers}. In the second case, $\mathcal{C}$ contains the middle term of the AR-sequence that starts in $(r_{j+1}, i_{j+1})$. In both cases, $\mathcal{C}$ is not submodule closed.
		
		So by Definition \ref{DefBijection}, there is some $\underline{w}$ with 
		\begin{equation*}
		\rho(\underline{w}) = (r_1, i_1)(r_2, i_2) \dots (r_m, i_m)
		\end{equation*}
		and $\mathcal{C} = \mathcal{C}_{\underline{w}}$.
	\end{proof}
	
	Recall that $\mathcal{C}_w = \mathcal{C}_{\underline{w}}$, where $\underline{w}$ is the leftmost word for $w$. So we need to prove that the word $\underline{w}$ in Lemma \ref{LemWord} is leftmost. Furthermore, we need the other direction, namely, that $\mathcal{C}_{\underline{w}}$ is submodule closed if $\underline{w}$ is leftmost.
	
	We will use the following lemma for the proofs of both directions:
	
	\begin{l1} \label{Lem2Ind}
		Suppose that the words $\underline{w}$ and $\underline{w}''$ are equivalent and there are pairs $(r, h), (s, i), (t+1,j)$ and series of pairs $\rho_1, \rho_2, \rho_4$ so that
		\begin{equation*}
		\rho(\underline{w}) = \rho_1 (r, h) \rho_2 (s, i) \rho_3,
		\end{equation*}
		\begin{equation*}
		\rho(\underline{w}'') = \rho_1 \rho_2 (s, i) (t+1, j)\rho_4
		\end{equation*}
		and either $\rho_3 = \rho_4$, or a pair $(q, g)$ is in $\rho_4$ if and only if $(q-1, g)$ is in $\rho_3$.
		
		Furthermore, suppose that the word $\underline{x}$ with $\rho(\underline{x}) = \rho_1 (r, h)\rho_2$ is reduced and $m_{ij} \ge 3$.
		
		If $M_0=(s, i), M_1 =(t, j), M_2, \dots, M_{m_{ij}-3} \notin \mathcal{C}_{\underline{w}}$, then there are sequences of modules as in Lemma \ref{LemExSeq} (used on $M_0$ and any $U \in \mathcal{C}_{\underline{w}}$) that yield some $U', Y \in \mathcal{I}$ so that
		\begin{equation} \label{ExactSequence}
		\xymatrix{0 \ar[r] &M_0 \ar[r] &M_{m_{ij}-2}\oplus U'\ar[r] & Y \ar[r] &0}
		\end{equation}
		is an exact sequence and either $Y \in \mathcal{C}_{\underline{w}}$ or both $Y = (r, h)^{E(m_{ij}-3)}$ and $U' \in \mathcal{C}_{\underline{w}}$ hold.
	\end{l1}
	
	\begin{proof}
		We show this by induction on the number $m$ of Coxeter relations needed to transform $\underline{w}$ into $\underline{w}''$.
		
		Furthermore, we show that a few additional assertions hold, which we need for the inductive proof:
		\begin{enumerate}[(\text{A}1)]
			\item If $(r', h')^{\gamma}\mid Y$ and $Z \in \mathcal{C}_{\underline{w}}$ is a direct summand of the middle term of the AR-sequence that ends in $(r', h')$, then $Z^{\gamma}\mid U'$.
			
			\item Let $V$ be the maximal direct summand of $U'$, so that for every $(q, g)\mid V$ there is some $(r', h')\mid Y$ and an indecomposable morphism $(q, g) \rightarrow (r', h')$. Furthermore, let $(q', g')$ be the biggest indecomposable direct summand of $V$ and $(r'', h'')$ the smallest indecomposable direct summand of $Y$.
			
			If there is some $o \in \mathbb{N}_0$ so that $Y, \tau^{-1} Y, \dots, \tau^{-o} Y \in \mathcal{C}_{\underline{w}}$ and $\tau^{-o-1} Y \notin \mathcal{C}_{\underline{w}}$, then one of the following holds: 
			\begin{enumerate}
				\item 	$\tau^{-o-1} Y = (r, h)^{E(m_{ij}-3)}$ and
				\begin{equation*}
				\tau^{-1} V, \tau^{-2} V, \dots, \tau^{-o-1} V \in \mathcal{C}_{\underline{w}}.
				\end{equation*} 
				
				\item Let $V'$ be the maximal direct summand of $X$ so that for every $(q, g) \mid V'$ there is some $0 \le k \le o$ so that $(q-k, g) \notin \mathcal{C}_{\underline{w}}$. Then there is some module $Y'$ with an exact sequence
				\begin{equation} \label{ExactSequence2}
				\xymatrix{0 \ar[r] &\tau^{-o-1} {V'} \ar[r] &\tau^{-o-1} Y\oplus U''\ar[r] &Y'\ar[r] &0}.
				\end{equation}
				Either $Y' \in \mathcal{C}_{\underline{w}}$ or both $Y' = (r, h)^{E(m_{ij}-3)}$ and $U'' \in \mathcal{C}'$ hold, where 
				\begin{equation*}
				\ind \mathcal{C}' = \ind \mathcal{C}_{\underline{w}} \setminus \mathcal{M}
				\end{equation*}
				with
				\begin{equation*}
				\mathcal{M} = \left\lbrace M \in \ind \mathcal{I} ~|~ \exists 0 \le o' \le o: \begin{Bmatrix} \tau^{o'+1} M \notin \mathcal{C}_{\underline{w}}\\
				(r'', h'') < \tau^{o'+1} M\\
				\tau^{o'+1} M < (q'-1, g')\\ 
				\end{Bmatrix}\right\rbrace 
				\end{equation*}
				Furthermore, (A1) and (A2) still hold if we exchange $U'$, $Y$ and $\mathcal{C}_{\underline{w}}$ for $U''$, $Y'$ and $\mathcal{C}'$ respectively.
			\end{enumerate}
		\end{enumerate}
		
		If there are some reflections $s_k, s_l$ and words $\underline{u}, \underline{v}$ so that
		\begin{equation*}
		w = \underline{u}\{s_ks_l\}^{m_{kl}}\underline{v} \quad \text{and} \quad w'' = \underline{u}\{s_ls_k\}^{m_{kl}}\underline{v},
		\end{equation*}
		then this is the result of Lemma \ref{LemRecSeq}. 
		
		Now suppose that $w = \underline{u}\{s_ks_l\}^{m_{kl}}\underline{v}$ and the lemma, (A1) and (A2) are proved for the word $\underline{w}_1= \underline{u}\{s_ls_k\}^{m_{kl}}\underline{v}$. 
		
		Either there are modules $U'_{(1)}, Y_{(1)}$ so that the exact sequence given by the inductive assumptions is
		\begin{equation}\label{ExactSequence3}
		\xymatrix{0 \ar[r] &M_0 \ar[r] &M_{m_{ij}-2}\oplus U'_{(1)}\ar[r] &Y_{(1)}\ar[r] &0},
		\end{equation}
		or we can write the exact sequence as a $\tau$-translate or a $\tau^{-1}$-translate of (\ref{ExactSequence3}).
		
		Let $\underline{w}_1 \equiv \underline{w}''_1$ with
		\begin{equation*}
		\rho(\underline{w}_1) = \rho'_1 (r_1, h_1) \rho'_2 (s', i) \rho'_3
		\end{equation*}
		and
		\begin{equation*}
		\rho(\underline{w}''_1) = \rho'_1 \rho'_2 (s', i) (t'+1, j)\rho'_4.
		\end{equation*}
		so that either $\rho'_3 = \rho'_4$, or a pair $(q, g)$ is in $\rho'_4$ if and only if $(q-1, g)$ is in $\rho'_3$.
		
		Furthermore, we can assume without loss of generality that $m_{kl}$ is even. Then there are $q_1, q_2 \in \mathbb{N}_0$ and a series of pairs $\rho''$ so that
		\begin{equation*}
		\rho(\underline{w}) = \rho(\underline{u}) (q_1- \frac{m_{kl}}{2}+1, k) (q_2 - \frac{m_{kl}}{2}+1, l) \dots (q_1 , k)(q_2, l) \rho''.
		\end{equation*}
		We can assume that $(r, h)$ is in the series of pairs $\rho(\underline{u}) (q_1- \frac{m_{kl}}{2}+1, k) (q_2 - \frac{m_{kl}}{2}+1, l) \dots (q_1 , k)(q_2, l)$, since otherwise there is nothing to show. Analogously, we assume that the pair $(r_1, h_1)$ is in the series of pairs $\rho(\underline{u}\{s_ls_k\}^{m_{kl}})$. Furthermore, we can assume $\underline{w}''_1 \neq \underline{w}''$ and $M_0 > (q_2, l)$
		
		Analogously to Lemma \ref{LemRecSeq}, if $m_{lk} \ge 3$, then there is some $X \in \mathcal{C}_{\underline{w}}$ so that
		\begin{equation} \label{ExactSequence4}
		\xymatrix{0 \ar[r] & (q_2, l) \ar[r] &  (q_1, k) \oplus X \ar[r] & (q_1- \frac{m_{kl}}{2}+1, k)}
		\end{equation}
		is an exact sequence.\\
		
		We have two different cases to consider:
		
		First, assume that $\underline{u}$ is also an initial subword of $\underline{w}''$. Then $\underline{u}\{s_ls_k\}^{m_{kl}-1}$ is an initial subword of $\underline{w}''_1$, since $\underline{x}$ with $\rho(\underline{x}) = \rho_1 (r, h) \rho_2$ is reduced. Furthermore, $\underline{w} <_l \underline{w}_1$ and $(q_1- \frac{m_{kl}}{2}+1, k) = (r, h)$. 
		
		If $(r_1, h_1) \in \mathcal{C}_{\underline{w}}$, then we can set $Y := Y_{(1)}$ and we have $(q_1, k) = (r_1 -1, h_1)$.
		
		(A1) holds by the inductive assumption, (A2) holds by (\ref{ExactSequence4}), (A1) and Lemma \ref{LemAlgInd}.
		
		On the other hand, if $(r_1, h_1) \notin \mathcal{C}_{\underline{w}}$, then $(q_1, k) = (r_1, h_1)$. Furthermore, if we have $Y_{(1)} = (r_1, h_1)^{E(m_{ij}-3)}$, then (A1), (\ref{ExactSequence4}) and Lemma \ref{LemAlgInd} give an exact sequence of the form (\ref{ExactSequence}) with $Y = (r, h)^{E(m_{ij}-3)}$. (A1) holds obviously. 
		
		If $Y_{(1)} \neq (r_1, h_1)^{E(m_{ij}-3)}$, then we set $Y := Y_{(1)}$ and we only need to prove that (A2) holds. Analogously to above, this is the result of Lemma \ref{LemAlgInd} and (\ref{ExactSequence4}).\\
		
		It remains to prove the assumption in the case that $\underline{u}$ is not an initial subword of $\underline{w}''$.
		
		If $\underline{w} <_l \underline{w}_1$, then $\rho(w_1)$ contains the series of pairs
		\begin{equation*}
		\rho(\underline{u}) (q_2 - \frac{m_{kl}}{2}+1, l) \dots (q_1 , k)(q_2, l)(q_1+1, k)
		\end{equation*}
		and the exact sequence given by the induction is either (\ref{ExactSequence3}) or the $\tau$-translate of (\ref{ExactSequence3}). We can assume without loss of generality, that some indecomposable direct summand of $Y_{(1)}$ or $\tau Y_{(1)}$ respectively is smaller than $(q_1+1, k)$. Otherwise, the arguments below hold analogously for an exact sequence given by (A2). 
		
		By Proposition \ref{Alg}, the exact sequence yielded by the induction is given by sequences of modules that fulfil (S1) - (S3). By Remark \ref{RemObs}, we can assume that $X''_{\gamma} = X'_{\gamma}$, for all $X'_{\gamma}$.
		
		In the following we begin with the case where this exact sequence is (\ref{ExactSequence3}). 
		
		By Lemma \ref{LemExSeq}, these series of modules yield an exact sequence
		\begin{equation} \label{ExactSequence5}
		\xymatrix{0 \ar[r] &M_0 \ar[r] & X_{\gamma} \oplus X'_{\gamma}\ar[r] & Y_{\gamma} \ar[r] &0}
		\end{equation}
		so that for every $(r'_1, h'_1) \mid Y_{\gamma}$ the inequality $(q_1, k) < (r'_1, h'_1)$ holds, since $M_0 > (q_1 + 1, k)$. We can even assume that $(q_1, k) < (r'_1, h'_1) \le (q_1+1, k)$ for all $(r'_1, h'_1) \notin \mathcal{C}_{\underline{w}_1}$. So either $(r'_1, h'_1) = (q_2, l)$ or $(r'_1, h'_1) = (q_1+1, k)$
		
		Furthermore, there is an irreducible morphism $X_{\gamma} \rightarrow (r'_1, h'_1)$ and $(q_1+1, k) \le X_{\gamma} \notin \mathcal{C}_{\underline{w}}$. Analogously to Lemma \ref{LemRecSeq}, if $m_{lk} \ge 3$, then there is some $X \in \mathcal{C}_{\underline{w}}$ so that
		\begin{equation}
		\xymatrix{0 \ar[r] & (q_1+1, k) \ar[r] &  (q_2, l) \oplus X \ar[r] & (q_1- \frac{m_{ij}}{2}+1, k)}
		\end{equation}
		with $V \in \mathcal{C}_{\underline{w}_1}$. Together with the $\tau$-translate of the exact sequence (\ref{ExactSequence4}), this shows that $X_{\gamma} \neq (q_2+1, l)$ and $X_{\gamma} \neq (q_1+1, k)$: Otherwise, by Lemma \ref{LemAlgInd}, we would get an exact sequence where either $(q_1- \frac{m_{kl}}{2}+1, k)$ or  $(q_2 - \frac{m_{kl}}{2}+1, l)$ is a direct summand of the last term, but every direct summand of the middle term is in $\mathcal{C}_{\underline{w}}$. This is a contradiction to the inductive assumption.		
		
		Since $(q_1 +2, k) > X_{\gamma} > (q_1+1, k)$, we have $\tau^{-1} X_{\gamma} \in \mathcal{C}_{\underline{w}_1}$. Inductively, $Y_{(1)} \in \mathcal{C}_{\underline{w}}$ and $(q_1, k) < (r'_1, h'_1)$ for every  $(r'_1, h'_1)\mid Y_{(1)}$.
		
		If $(q_1+1, k)\mid Y_{\gamma}$ for any exact sequence of the form (\ref{ExactSequence5}), then there is such a sequence so that $Y_{\gamma} \in \mathcal{C}_{\underline{w}}$, but $Y_{\gamma} \notin \mathcal{C}_{\underline{w}_1}$.
		
		Otherwise, the sequence (\ref{ExactSequence3}) is already of the form (\ref{ExactSequence}). 
		
		In the latter case, it is easily seen that this sequence fulfils (A1) and (A2): the former holds by the inductive assumption. Define $V$, $V'$, $\mathcal{C}'$ as in (A2) and let $V_{(1)}, V'_{(1)}$ be the corresponding modules, $\mathcal{C}'_{(1)}$ the corresponding category for the sequence (\ref{ExactSequence3}). Then $V = V_{(1)}$, $V' = V'_{(1)}$ and $\mathcal{C}' = \mathcal{C}'_{(1)}$. Thus, assertion (2) also holds.
		
		So assume that there is some $\gamma$ with $(q_1+1, k) \mid Y_{\gamma} \in \mathcal{C}_{\underline{w}}$. This sequence is of the form (\ref{ExactSequence}) and (A1) holds. Let $\alpha \in \mathbb{N}$ be the maximal exponent so that $(q_1+1, k)^{\alpha} \mid Y_{\gamma}$
		
		We can write $X_{\gamma} \oplus X'_{\gamma} = B_1 \oplus B'_1 \oplus M_{m_{ij}-2}$ so that $B_1$ is the maximal direct summand of $X_{\gamma} \oplus X'_{\gamma}$ with an irreducible morphism $B_1 \rightarrow (q_1+1, k)$. Then there is an exact sequence
		\begin{equation*}
		\xymatrix{0 \ar[r] & B_1 \ar[r] & C_1 \oplus  (q_1+1, k)^{\alpha} \ar[r] & \tau^{-1} B_1}.
		\end{equation*}
		By Remark \ref{RemObs} and Lemma \ref{LemAlgInd}, $U'_{(1)} = B'_1 \oplus C_1$. If $Y_{\gamma} = D \oplus (q_1+1, k)^{\alpha}$, then $Y_{(1)} = D \oplus \tau^{-1} B_1$.
		
		So we can write $X'_{(1)} = B''_1 \oplus C''_1$ with $B''_1\mid B'_1$ and $C''_1\mid C_1$. We get $X' = B''_1 \oplus B_1$. By Proposition \ref{Alg}, Lemma \ref{LemAlgInd} and the inductive assumption, assertion (2) is fulfilled.\\
		
		If we still have $\underline{w} <_l \underline{w}_1$, but the exact sequence given by the inductive assumption is the $\tau$-translate of (\ref{ExactSequence3}), then analogously we have $\tau Y_{(1)} \in \mathcal{C}_{\underline{w}_1}$ and $(q_1, k) < (r'_1, h'_1)$ for all direct summands $(r'_1, h'_1)$ of $\tau Y_{(1)}$. Suppose that we have $\tau Y_{(1)}, Y_{(1)}, \dots, \tau^{-o+1} Y_{(1)} \in \mathcal{C}_{\underline{w}}$ and $\tau^{-o} Y \notin \mathcal{C}_{\underline{w}}$. If $o > 0$, then (\ref{ExactSequence3}) is of the form (\ref{ExactSequence}). By the inductive hypothesis, (A2) is fulfilled.
		
		If $o = 0$, then we use assertion (A2) of the inductive hypothesis: by Lemma \ref{LemAlgInd}, there is an exact sequence of the form (\ref{ExactSequence}) that fulfils (A1) and (A2).\\
		
		It only remains to regard what happens if $\underline{w}_1 <_l \underline{w}$. Then $\rho(w_1)$ contains the series of pairs
		\begin{equation*}
		\rho(\underline{u}) (q_2 - \frac{m_{kl}}{2}, l) \dots (q_2-1, l)(q_1 , k).
		\end{equation*}
		The exact sequence given by the inductive assumption is either (\ref{ExactSequence3}) or the $\tau^{-1}$-translate. In the first case, we show analogously to the above that $Y_{(1)} \in \mathcal{C}_{\underline{w}_1}$ and $(q_1, k) < (r'_1, h'_1)$ for all $(r'_1, h'_1) \mid Y_{(1)}$. If $(q_2, l)\nmid Y_{(1)}$, then the sequence (\ref{ExactSequence3}) is already of the form (\ref{ExactSequence}) and the assertions hold.
		
		Otherwise, $\tau^{-1} Y_{(1)} \notin \mathcal{C}_{\underline{w}_1}$. Analogously to before, $(q_1, k) \nmid \tau^{-1} X'_{(1)}$ by Lemma \ref{LemRecSeq} and (A2): If $(q_1, k) \mid \tau^{-1} X'_{(1)}$ we use Lemma \ref{LemAlgInd} and get an exact sequence where either $(q_1- \frac{m_{kl}}{2}, k)$ or  $(q_2 - \frac{m_{kl}}{2}, l)$ is a direct summand of the last term, but every direct summand of the middle term is in $\mathcal{C}_{\underline{w}}$. This is a contradiction to the inductive assumption.
		
		So $\tau X'_{(1)} \in \mathcal{C}_{\underline{w}}$, which means that no direct summand of $X'_{(1)}$ is in $\mathcal{C}_{\underline{w}}$. 
		
		As before, the assertion (A2) of the inductive hypothesis and Lemma \ref{LemAlgInd} show that there is an exact sequence of the form (\ref{ExactSequence}) and that (A1) and (A2) are fulfilled.\\
		
		If the exact sequence yielded by the inductive hypothesis is the $\tau^{-1}$-translate of (\ref{ExactSequence3}), then similar to the arguments above we get an exact sequence (\ref{ExactSequence}) so that $Y \in \mathcal{C}_{\underline{w}}$ and $(q_1, k) < (r', h')$ for every direct summand $(r', h')$ of $Y$. This exact sequence fulfils (A1) and (A2).
	\end{proof}

	\section{The first direction} \label{SecProof2}
	
	In this section we show inductively that for every $w \in W$, the category $\mathcal{C}_{w}$  is submodule closed. Afterwards, it only remains to show that every cofinite, submodule closed category is of the form $\mathcal{C}_w$.
	
	We begin with the basis of the induction:
	
	\begin{l1} \label{LemStartInd}
		Let $m_{ij} < \infty$ and $U_1, \dots, U_{m_{ij}-1}$ be as in Lemma \ref{LemRecSeq}. If we have $U_1, \dots, U_{m_{ij}-1} \in \mathcal{C}_{\underline{w}}$ and  $M_0 \notin \mathcal{C}_{\underline{w}}$, then $\mathcal{C}_{\underline{w}}$ is not submodule closed and $\underline{w}$ is not leftmost. 
	\end{l1}
	
	\begin{proof}
		Since $M_0 \subseteq U_{m_{ij}-1}$, the category $\mathcal{C}_{\underline{w}}$ is not submodule closed. Let $\underline{w}$ be a word for the element $w \in W$. By Definition \ref{DefMs}, $M_0=(s, i)$ and $M_1=(t,j)$. 
		
		Suppose that $M_1 \in \mathcal{C}_{\underline{w}}$. Let $\underline{u}$ be the initial subword of $\underline{w}$ that is defined through the inequality $(r,k) \le (s-1,i)$ for every pair $(r,k)$ in  $\rho(\underline{u})$.
		
		Then there are reflections $s_{k_1}, s_{k_2}, \dots, s_{k_m}$ and a word $\underline{v}$ so that
		\begin{equation*}
		\underline{w} = \underline{u}s_{k_1} s_{k_2} \dots, s_{k_m} s_i\underline{v}.
		\end{equation*}
		Since $U_1 \in \mathcal{C}_{\underline{w}}$, we have
		\begin{equation*}
		m_{k_1 i} = m_{k_2 i} = \dots = m_{k_m i} = 2
		\end{equation*}
		by Lemmas \ref{LemARWeyl} and \ref{LemARNumbers}. So
		\begin{equation*}
		\underline{w'} = \underline{u}s_i s_{k_1} s_{k_2} \dots s_{k_m} \underline{v}
		\end{equation*} 
		is equivalent to $\underline{w}$ and thus a word for $w$.
		Since $(r,k) \le (s-1,i)$ for all reflections $(r,k)$ in $\rho(\underline{u})$, we see that either $\underline{w}' <_l \underline{w}$, or $\underline{w}$ is not reduced.
		
		Clearly, the same argument holds if $M_m \in \mathcal{C}_{\underline{w}}$ for some $1 \le m \le m_{ij} - 1$.
		
		It remains to prove that $\mathcal{C}_{\underline{w}}$ is not submodule closed if $M_0, \dots, M_{m_{ij}-1} \notin \mathcal{C}_{\underline{w}}$. Without loss of generality, we can assume that $M_{m_{ij}-1} = (p,i)$ for some $p \in \mathbb{N}$ and $M_{m_{ij}-2} = (q,j)$ for some $q \in \mathbb{N}$. Suppose that $(r,k)$ is a pair in $\underline{w}$. If $(q-1,j) < (r,k) < (t,j)$ then we use that $U_1, U_2, \dots U_{m_{ij}-1} \in \mathcal{C}_{\underline{w}}$ and get $m_{jk} = 2$ by \ref{CorRecSeq}. If $(p,i) < (r,k) < (s, i)$ then $m_{ik} = 2$.
		
		Let $\underline{u}'$ be the initial subword of $\underline{w}$ that is defined through the inequality $(r,k) \le (q-1,j)$ for every pair $(r,k)$ in  $\rho(\underline{u})$.
		
		Then there are reflections $s_{k_1}, \dots s_{k_m}$ with
		\begin{equation*}
		m_{k_1 j} = m_{k_2 j} = \dots = m_{k_m j} = 2
		\end{equation*}	
		so that $\underline{w} \equiv \underline{w}'$ for
		\begin{equation}
		\underline{w}' = \underline{u}' s_{k_1} \dots s_{k_m}\{s_is_j\}^{m_{ij}} \underline{v}'.
		\end{equation}
		
		So $\underline{w}$ is also equivalent to
		\begin{equation}
		\underline{w}'' = \underline{u}' s_j  s_{k_1} \dots s_{k_m}\{s_is_j\}^{m_{ij}-1} \underline{v}'
		\end{equation}
		and either $\underline{w}'' <_l \underline{w}$ or $\underline{w}$ is not reduced.
	\end{proof}
	
	This proof even shows the following:
	
	\begin{c1} \label{CorStartInd}
		If $m_{ij} < \infty$, $M_0, M_1, \dots, M_{m_{ij}-1} \notin \mathcal{C}_{\underline{w}}$ and $\underline{w}$ is reduced, then there is a word $\underline{w}' \equiv \underline{w}$ with $\underline{w}' <_l \underline{w}$, pairs $(r, h) = M_{m_{ij}-1}$, $(r', h') \neq (r, h)$ and series of pairs $\rho_1, \rho_2, \rho_3$ so that 
		\begin{equation*}
		\rho(\underline{w}') = \rho_1 (r, h) \rho_2 \quad \text{and} \quad 
		\rho(\underline{w}) = \rho_1 (r', h')\rho_3.
		\end{equation*}
	\end{c1}
	
	For the inductive step, we still need some lemmas:
	
	\begin{l1} \label{LemMijinWord}
		Suppose that for some $\underline{w}$, we have $M_0 \notin \mathcal{C}_{\underline{w}}$ and  $M_0$ is a submodule of $U \in \mathcal{C}$. Let $U_{m_{ij}-1}$ be as in Lemma \ref{LemRecSeq} with modules $(r_k, l_k) \notin \mathcal{C}_{\underline{w}}$ for $1 \le k \le a$ so that $\bigoplus_{k=1}^a (r_k, l_k) \mid U_{m_{ij}-1}$.
		
		Let 
		\begin{align} \label{Modulseq}
		\begin{split}
		&(X_1, X_2, \dots, X_m)\\
		&(X'_1, X'_2,\dots, X'_m)\\
		&(Y_1, Y_2, \dots, Y_m)
		\end{split}
		\end{align}
		be the sequences of modules that yield the exact sequences
		\begin{equation*}
		\eta_k: \xymatrix{0 \ar[r] & M_0 \ar^-{f_{k}}[r] & M_{k}^{E(k)} \oplus U_{k} \ar[r] & M_{k+1}^{E(k-1)} \ar[r] & 0}
		\end{equation*}
		for all $ 1 \le k \le m_{ij}-1$. 
		
		Then one of the following holds:
		\begin{enumerate}[(a)]
			\item There is some $1 \le m' \le m_{ij}-1$  and some $U' \in \mathcal{C}_{\underline{w}}$ with a monomorphism $M_{m'} \rightarrowtail U'$
			
			\item If $(X_1, \dots X_m, X_{m+1}, \dots X_{m'})$ is part of a triple of sequences that fulfils (S1) - (S5), there is some $1 \le k \le m'$ so that $M_{m_{ij}} \mid X_{k} \oplus X'_{k}$. Furthermore, for $l \in \mathbb{N}$, $1 \le k \le a$ and $M_{m_{ij}-1} < (r_k-l, l_k)$, we have $(r_k -l, l_k) \notin \mathcal{C}_{\underline{w}}$. 
			
			\item $a = 1$, $m_{il_1}+m_{jl_1} = 5$ and there is no indecomposable morphism $M_{m''} \rightarrow (r_1, l_1)$ for $m'' < m_{ij}-3$.
			
			\item $a = 1$ and $(r_1, l_1) < M_{m_{ij}-1}$.
		\end{enumerate}
		If  $(r_k-1, l_k) \in \mathcal{C}_{\underline{w}}$ for some $1 \le k \le a$, then (a) holds.
	\end{l1}
	
	\begin{proof}
		By Corollary \ref{CorRecSeq}, for all $1 \le k \le a$, there are some $X \in \mathcal{I}$, $\beta_k \in \mathbb{N}$ and $2 \le m_k \le m_{ij}$ so that the AR-sequence that starts in $(r_k, l_k)$ is of the form
		\begin{equation*}
		\xymatrix{0 \ar[r]& (r_k, l_k) \ar[r] & M^{\beta_k}_{m_k} \oplus X \ar[r]& (r_k -1, l_k) \ar[r] & 0}.
		\end{equation*}

		First, suppose that $a > 1$.	 By Lemma \ref{LemRec} and Lemma \ref{LemRecSeq}, there is an exact sequence
		\begin{equation*}
		\xymatrix{0 \ar[r] & M_0 \ar[r] & U_{m_{ij}-1} \ar[r] & M_{m_{ij}} \ar[r] & 0}.
		\end{equation*}
		and we can assume that $X_m \oplus X'_m =  U_{m_{ij}-1}$ and $Y_m = M_{m_{ij}}$. By Corollary \ref{CorAlg}, we can set $X_m := (r_1, l_1)$ and $X_{m+1} := (r_2, l_2)$. If $M_{m_1} = M_{m_2} = M_{m_{ij}}$, then $M_{m_{ij}} \mid X_{m+2} \oplus X'_{m+2}$ by (S3). Thus (b) is fulfilled.
		
		Otherwise, set $\gamma = m+m_{ij}+2-m_1$ , $\delta = \gamma + m_{ij} - m_2$
		\begin{equation*}
		X_{m+3} := M_{m_1}, X_{m+4} := M_{m_1+ 1}, \dots, X_{\gamma} := M_{m_{ij}-1}
		\end{equation*}
		and
		\begin{equation*}
		X_{\gamma+1} := M_{m_2}, X_{\gamma + 2} := M_{m_2+ 1}, \dots X_{\delta} := M_{m_{ij}-1}.
		\end{equation*}
		Then $M_{m_{ij}} \mid X_{\delta+1} \oplus X'_{\delta+1}$ and (b) holds.
		
		On the other hand, suppose that there is an indecomposable morphism $M_{m''} \rightarrow (r_1, l_1)$ for some $m'' < m_{ij}-3$. By Corollary \ref{CorRecSeq}, we have $(r_1-1, l_1) \mid U_{m_{ij}-1}$. By the construction of $\eta_{m''+1}$ from $\eta_{m''+1}$, there must be a module $U'_{m''+1}$ so that we can write $\eta_{m''+1}$ as the following:
		\begin{equation*}
		\xymatrix{0 \ar[r] & M_0 \ar[r] & M_{m''+1}^{E(m''+1)} \oplus (r_1, l_1)^{E(m'')} \oplus U'_{m''+1}\ar[r] &M_{m''+2}^{E(m'')}}.
		\end{equation*}
		So there is a module $U'$ with an exact sequence
		\begin{equation*}
		\xymatrix{0 \ar[r] & M_0 \ar[r] &  M_{m''+1}^{E(m''+1)} \oplus U'\ar[r] &(r_1-1, l_1)^{E(m''+1)} \ar[r]& 0}.
		\end{equation*}
		By Corollary \ref{CorMono} and the monomorphism $M_0 \rightarrowtail U$, we get a monomorphism 
		\begin{equation*}
		M_{m_{ij}-1}\rightarrowtail U \oplus (r_1-1, l_1).
		\end{equation*}
		So either (a) is fulfilled or $(r-1, l_1) \notin \mathcal{C}_{\underline{w}}$. In this case, $a > 1$ and (b) is fulfilled.
		
		Finally, suppose that $m_{il_1} + m_{jl_1} > 5$ and $a = 1$.  By Lemma \ref{LemARWeyl}, $\beta_1 > 1$ if either $m_{il_1} \ge 4$ or $m_{il_2} \ge 4$. 	
		We can set $X_{m+1} := (r_1, l_1)$ and $X_{m+2} := X_{m+3} := M_{m_1}$. Analogously to the case $a > 1$, this yields some $m' \in \mathbb{N}$ so that $M_{m_{ij}} \mid X_{m'} \oplus X'_{m'}$.
		
		If $m_{il_1} = m_{jl_1} = 3$, then either (d) holds, or we can write the AR-sequence that starts in $(r_1, l_1)$ as
		\begin{equation} \label{ARSeq}
		\xymatrix{0 \ar[r]& (r_1, l_1) \ar[r] & M_{m_1} \oplus M_{m_1+1} \oplus X' \ar[r]& (r_1 -1, l_1) \ar[r] & 0 }
		\end{equation}
		for some $X' \in \mathcal{I}$ and some $1 \le m_1 \le m_{ij}-1$. Since $X_{m+1} := (r_1, l_1)$, $X_{m+2} := M_{m_1}$ and $X_{m+3} = M_{m_1+1}$, we see that (b) holds.
	\end{proof}
	
	We can generalize this lemma in the following way, which we will need for an induction:
	
	\begin{r1} \label{RemMijInWord}
		Let us assume that the assumptions of Lemma \ref{LemMijinWord} hold, except that $M_0$ is not a submodule of some $U \in \mathcal{C}_{\underline{w}}$ but of $X \oplus U$ for some indecomposable $X \notin \mathcal{C}_{\underline{w}}$. If $X = M_k$ for some $1 \le k \le m_{ij}-1$, we alter the assumptions on the sequences (\ref{Modulseq}) accordingly so that they fulfil (S1) - (S3) with respect to $M_0$ and $X \oplus U$. 
		
		Furthermore, we assume that $X^2 \mid X_k \oplus X'_k$ for some $1 \le k \le m$: then for every sequence $(X_1, \dots, X_m, X_{m+1}, \dots, X_{m'})$ which is part of a triple of sequences that fulfils (S1) - (S5), we have some $1 \le k \le m'$ so that $X_k = X$. 
		
		So every argument in the proof of \ref{LemMijinWord} still holds and we still get that one of the conditions (a) - (d) must be fulfilled.
		
		If we take a look at case (b) of \ref{LemMijinWord} and suppose that $a = 1$, then by Lemma \ref{LemAlgInd}, there is a monomorphism
		$(r_1, l_1) \rightarrowtail M_{m_{ij}} \oplus U$. So there	are sequences of modules
		\begin{align*}
		\begin{split}
		&(^1X_1, ^1X_2, \dots, ^1X_{o})\\
		&(^1X'_1, ^1X'_2,\dots, ^1X'_{o})\\
		&(^1Y_1, ^1Y_2, \dots, ^1Y_{o})
		\end{split}
		\end{align*}
		which fulfil (S1) - (S5) with respect to $(r_1, l_1)$ and $M_{m_{ij}} \oplus U$. Since $M_{m_{ij}} \mid X_k \oplus X'_k$ for some $1 \le k \le m'$, there must be some $k'$ so that $M^2_{m_{ij}} \mid ^1X_{k'} \oplus X'_{k'}$.
	\end{r1}
	
	We still need a result about the case $a > 1$:
	
	\begin{l1} \label{LemA=2}
		Suppose that for some $\underline{w}$, we have $M_0 \notin \mathcal{C}_{\underline{w}}$ and  $M_0$ is a submodule of $U \in \mathcal{C}$. Let $U_{m_{ij}-1}$ be as in Lemma \ref{LemRecSeq} with modules $(r_k, l_k) \notin \mathcal{C}_{\underline{w}}$ for $1 \le k \le a$ so that $\bigoplus_{k=1}^a (r_k, l_k) \mid U_{m_{ij}-1}$.
		
		Furthermore, suppose that for all $1 \le k \le a$, there is an irreducible morphism $M_m \rightarrow (r_k, l_k)$ for some $0 \le m < m_{ij}-2$.
		
		If $a > 1$, then one of the following holds:
		
		\begin{enumerate}[(a)]
			\item There is some $N < M$, $N \notin \mathcal{C}_{\underline{w}}$ that is a submodule of some $U' \in \mathcal{C}_{\underline{w}}$. 
			
			\item We have $a = 2$,  $m_{l_1l_2} = 2$, $m_{il_k} = 3$, $m_{jl_k} = 2$ for $1 \le k \le 2$ and $m_{ij} = 3$.
			
			\item We have $a = 2$, $l_1 = l_2$, $m_{il_1} = 3$ and $m_{ij} = 3$.
		\end{enumerate}
		In case (b),
		\begin{equation*}
		M_{m_{ij}}, M_{m_{ij}+1}, M_{m_{ij}+2}, M_{m_{ij}+3} \notin \mathcal{C}_{\underline{w}}.
		\end{equation*}
		and
		\begin{equation*}
		(r_1-1, l_1), (r_2-1, l_2),  (r_1-2, l_1), (r_2-2, l_2) \notin \mathcal{C}_{\underline{w}}.
		\end{equation*}
	\end{l1}
	
	\begin{proof}	
		The proof is analogous to that of Lemma \ref{LemRecSeq}. First note that if $m_{ij} = 3$, we get $m_{il_k} \ge 3$ for $1 \le k \le a$, since there is an irreducible morphism $M_m \rightarrow (r_k, l_k)$ for some $0 \le m < m_{ij}-2$. 
		
		If $m_{l_1l_2} \neq 2$, we can exchange $j$ and $l_2$ in the calculations below.	
		
		Define $\alpha_{ki}$, $\beta_{ki}$, $\alpha_{kj}$, $\beta_{kj}$ analogous to $\alpha$ and $\beta$ with $i, l_k$ and $j, l_k$ instead of $i, j$. If $m_{l_1l_2} = 2$, but (b) is not fulfilled, then
		\begin{equation*}
		\alpha \beta + \sum_{k=1}^a (\alpha_{ki}+\alpha_{kj})(\beta_{ki}+\beta_{kj}) \ge 4.
		\end{equation*}
		By Corollary \ref{CorRemAlg}, we can do completely analogous calculations to the case $\alpha \beta = 4$. In these calculations, we construct exact sequences which contain modules of the form $M_o$, $(o', l_1)$ and $(o'', l_2)$ for some $o, o', o'' \in \mathbb{N}$.

		It remains to show that (a) is fulfilled if any of these modules is in $\mathcal{C}_{\underline{w}}$. 
		
		By Lemma \ref{LemMijinWord}, if $(r_k-1, l_k) \in \mathcal{C}_{\underline{w}}$ for some $1 \le k \le a$, then (a) holds. The rest follows inductively with the same argument as in \ref{LemMijinWord}.
	\end{proof}
	
	Now we can show the following, which is the last lemma that we need to prove the first direction of the main theorem:
	\begin{l1} \label{LemWordMij}
		Let $\underline{w}$ be a word so that $(s, i) = M_0, M_1, \dots, M_{m_{ij}-1} \notin \mathcal{C}_{\underline{w}}$. Then there are words $\underline{u}, \underline{v}$ so that $\underline{w} = \underline{u} s_i \underline{v}$ and there is some $\rho$ with $\rho(\underline{w}) = \rho(\underline{u}) (s, i) \rho$. Suppose that there is some $U \in \mathcal{C}_{\underline{w}}$ with a monomorphism $M_0 \rightarrowtail U$ and for every $X < M_0$ with some $U' \in \mathcal{C}_{\underline{w}}$ and a monomorphism $X \rightarrowtail U'$, we have $X \in \mathcal{C}_{\underline{w}}$.
		
		Then there exists some $\underline{u}'$ so that
		\begin{equation} \label{EqWordChangeable}
		\underline{w} \equiv \underline{u}'\{s_is_j\}^{m_{ij}} \underline{v}.
		\end{equation} 
	\end{l1}
	
	\begin{proof}	
		By Lemma \ref{LemRec} and \ref{LemRecSeq}, Proposition \ref{Alg} yields an exact sequence 
		\begin{equation} \label{SeqUMij}
		\xymatrix{0 \ar[r] & M_0 \ar[r] & U_{m_{ij}-1} \ar[r] & M_{m_{ij}} \ar[r] & 0}.
		\end{equation}
		If $(r, l) \in \mathcal{C}_{\underline{w}}$ for all $(r, l) \mid U_{m_{ij}-1}$ with $(r, l) > M_{m_{ij}-1}$, then (\ref{SeqUMij}) is obvious by Corollary \ref{CorRecSeq}.
		
		Otherwise, we get $M_{m_{ij}} \notin \mathcal{C}_{\underline{w}}$, since there is a monomorphism $U_{m_{ij}-1} \rightarrowtail M_{m_{ij}} \oplus U$ by Lemma \ref{CorMono}.

		There are direct summands $(r_1, l_1), \dots, (r_a, l_a) \notin \mathcal{C}_{\underline{w}}$ of $U_{m_{ij}-1}$ so that $M_{m_{ij}-1} < (r_k, l_k)$ for all $1 \le k \le a$. 
		
		We can assume without loss of generality that $m_{ij}$ is odd; otherwise we only need to exchange $s_j$ and $s_i$ in the arguments below.
		
		If there is no morphism  $M_m \rightarrow (r_k, l_m)$ for some $0 \le m < m_{ij} - 2$ for all $1 \le k \le a$, then $m_{l_ki} = 2$ and for some word $\underline{x}$ we have
		\begin{equation*}
		\underline{w} = \underline{x} s_i s_{l_1} \dots s_{l_a} \{s_js_i\}^{m_{ij}-1} \equiv \underline{x} s_{l_1} \dots s_{l_a} \{s_is_j\}^{m_{ij}-1}.
		\end{equation*}
		
		So we can assume that for $(r_k, l_k)$, there is a morphism $M_m \rightarrow (r_k, l_m)$ for some $0 \le m < m_{ij} - 2$.
		
		By Lemma \ref{LemMijinWord}, one of the following cases hold:
		\begin{enumerate}[(a)]
			\item $a = 1$, $m_{il_1}+m_{jl_1} = 5$ and there is no indecomposable morphism $M_{m''} \rightarrow (r_1, l_1)$ for $m'' < m_{ij}-3$.
			
			\item There is some $m'$ so that $M_{m_{ij}} \mid X_{m'} \oplus X'_{m'}$
		\end{enumerate}	
		In case (a) there are words $\underline{u}'', \underline{v}$, so that either
		\begin{equation} \label{Form1}
		\underline{w} = \underline{u}'' s_is_{l_1}\{s_js_i\}^{m_{ij}-1}\underline{v}
		\end{equation}
		or both $m_{jl_1} = 2$ and 
		\begin{equation} \label{Form2}
		\underline{w}=\underline{u}''s_is_js_{l_1}\{s_is_j\}^{m_{ij}-2}\underline{v}.
		\end{equation}
		Then there is some word $\underline{u}_1$ so that $\underline{w} \equiv \underline{u}_1 s_{l_1} \{s_is_j\}^{m_{ij}}\underline{v}$: If $m_{il_1} = 2$, there is nothing to show. If $m_{jl_1} = 2$, then we use that $M_{m_{ij}}$ is of the form $(q, j)$ for some $q$ and there is monomorphism $U_{m_{ij}-1} \rightarrowtail M_{m_{ij}} \oplus U$:	Let $\underline{u}_2$ be the subword of $\underline{u}_1$ that does not contain the reflection that corresponds to $(q, j)$. By the inductive assumption, $\underline{u}_2 s_i s_{l_1}$ is equivalent to a smaller word, since $\mathcal{C}_{\underline{u}_2s_is_{l_1}}$; thus it is equivalent to a word $\underline{u}_3\{s_is_{l_1}\}^{m_{il_1}}$.
		
		Now we go back to looking at $\underline{u}_1$, not $\underline{u}_2$. Since $m_{il_1} = 3$ and $m_{jl_1}$, there is some word $\underline{u}_4$ with $\underline{u}_1s_is_{l_1} \equiv \underline{u}_4 \{s_is_{l_1}\}^{m_{il_1}}$ and thus
		\begin{equation*}
		\underline{w} \equiv \underline{u}'' s_{l_1} \{s_is_j\}^{m_{ij}}\underline{v}
		\end{equation*}
		for some word $\underline{u}''$.
		
		On the other hand, suppose that 
		\begin{equation} \label{EqMij}
		M_{m_{ij}} \mid X_{m'} \oplus X'_{m'}.
		\end{equation}
		Since we have $\bigoplus_{k=1}^{a}(r_k, l_k) \mid U_{m_{ij}}$, there is a monomorphism
		\begin{equation} \label{EqMono}
		\bigoplus_{k=1}^{a}(r_k, l_k) \rightarrowtail U \oplus M_{m_{ij}}.
		\end{equation}
		Assume that $(r_1, l_1) \le (r_2, l_2) \le \dots \le (r_a, l_a)$. 
		
		Let the AR-sequence that starts in $(r_a, l_a)$ be
		\begin{equation*}
		\xymatrix{0 \ar[r] & (r_a, l_a) \ar[r] & M_{m_1} \oplus Z_1 \ar[r] & (r_a-1, l_1)\ar[r] &0}.
		\end{equation*}
		
		If $a = 1$, then by Lemma \ref{LemMijinWord}, $\underline{w}$ must have the form (\ref{Form1}) or (\ref{Form2}). The only difference to case (a) is that $m_{il_1}+m_{jl_2} > 5$.
		
		We can use case (a) as the basis of an induction: Instead of the modules $M_0, M_1, (r_1, l_1)$, we take $(r_1, l_1), M_{m_1}, M_{m_1+1}$ and use the same arguments as before. Since $M_{m_{ij}} \mid X_{m'} \oplus X'_{m'}$, we can use Remark \ref{RemMijInWord} and either we get the analogue to case (a) above or the analogue to the case (b). In the first case, we get $m_{il_1} + m_{ij} = 5$ or $m_{jl_1} + m_{ij} = 5$ and there is some $\underline{u}'_1$ so that $\underline{w}$ is equivalent to a word with the subword $\{s_{l_1}s_{i}\}^{m_{il_1}}$ if $m_1$ is even and $\{s_{l_1}s_{j}\}^{m_{jl_1}}$ is odd. Thus we also get $\underline{w} \equiv \underline{u}' \{s_is_j\}^{m_{ij}} \underline{v}$ for some words $\underline{u}', \underline{v}$. In the case (b), we continue this inductively.
		
		After finitely many steps we get
		\begin{equation*}
		\underline{w} \equiv \underline{u}' \{s_is_j\}^{m_{ij}} \underline{v}.
		\end{equation*}

		If $a \neq 1$, then $a = 2$ by Lemma \ref{LemA=2}. If $l_1 = l_2$, then $l_1 = l_2$, $m_{il_1} = 3$ and $m_{ij} = 3$. We can exchange $j$ and $l_1$ to get the case $a= 1$. Otherwise, $m_{l_1l_2} = 2$, $m_{il_k} = 3$, $m_{jl_k} = 2$ for $1 \le k \le 2$ and $m_{ij} = 3$. Furthermore,
		\begin{equation*}
		M_{m_{ij}}, M_{m_{ij}+1}, M_{m_{ij}+2}, M_{m_{ij}+3} \notin \mathcal{C}_{\underline{w}}.
		\end{equation*}
		and
		\begin{equation*}
		(r_1-1, l_1), (r_2-2, l_2),  (r_1-2, l_1), (r_2-2, l_2) \notin \mathcal{C}_{\underline{w}}.
		\end{equation*}
		Analogously to before, we see inductively that $\underline{w}$ is equivalent to a word with the subword 
		\begin{equation*}
		s_{l_1}s_{l_2}s_js_is_{l_1}s_{l_2}s_js_is_{l_1}s_{l_2}s_js_i.
		\end{equation*}
		(For the purpose of this induction, we can treat the word above completely analogously to a word of the form $\{s_is_j\}^{m_{ij}}$ with $m_{ij} = 6$. As in \ref{LemA=2}, all calculations are the same by Corollary \ref{CorRemAlg}.)
		
		We have the following equivalences, where bold reflections denote those which differ from the reflections in the word above:	
		\begin{align*}
		&s_{l_1}s_{l_2}s_js_is_{l_1}s_{l_2}s_js_is_{l_1}s_{l_2}s_js_i\\
		\equiv &s_{l_1}s_{l_2}s_js_i\mathbf{s_js_{l_1}s_{l_2}}s_i\mathbf{s_{l_2}s_{l_1}}s_js_i\\
		\equiv &s_{l_1}s_{l_2}\mathbf{s_is_js_i}s_{l_1}\mathbf{s_{i}s_{l_2}s_{i}}s_{l_1}s_js_i\\
		\equiv &s_{l_1}s_{l_2}s_is_j\mathbf{s_{l_1}s_{i}s_{l_1}}s_{l_2}s_{i}s_{l_1}s_js_i\\
		\equiv &s_{l_1}s_{l_2}s_is_js_{l_1}s_{i}\mathbf{s_{l_2}s_{l_1}}s_{i}s_{l_1}s_js_i\\
		\equiv & s_{l_1}s_{l_2}s_is_js_{l_1}s_{i}s_{l_2}\mathbf{s_{i}s_{l_1}s_{i}}s_js_i.
		\end{align*}
		So $\underline{w}$ is equivalent to a word with the subword $s_is_js_i$ and the assertion is true.
	\end{proof}
	
	Finally, we can prove the first direction of our main result:
	
	\begin{l1} \label{LemFirstDirection}
		If $\underline{w}$ is a leftmost word, then $\mathcal{C}_{\underline{w}}$ is submodule closed.
	\end{l1}
	
	\begin{proof}
		Suppose that $\mathcal{C}_{\underline{w}}$ is not submodule closed. Then there is some $M_0 \in \ind \mathcal{I} \setminus \mathcal{C}_{\underline{w}}$ and some $U \in \mathcal{C}_{\underline{w}}$ with a monomorphism $M_0 \rightarrowtail U$. Furthermore, we can assume that for every $X < M_0$ with some $U' \in \mathcal{C}_{\underline{w}}$ and a monomorphism $X \rightarrowtail U'$, we have $X \in \mathcal{C}_{\underline{w}}$.
		
		We use induction on the length $m$ of the sequences of modules in Proposition \ref{Alg} applied on $M_0$ and $U$. If $m =1$, then $\underline{w}$ is not leftmost by Lemma \ref{LemStartInd}. 
		
		Now suppose that $\underline{w}$ is not leftmost if the sequences have the length $m$ or smaller. We prove that this is also the case if they have length $m+1$:
		
		We can assume without loss of generality that $M_1 \notin \mathcal{C}_{\underline{w}}$, since $m +1 > 1$. On the other hand, by Lemma \ref{LemRecSeq}, the sequences of modules induce an exact sequence
		\begin{equation*}
		\xymatrix{0 \ar[r] & M_0 \ar[r] & M_1^{E(1)} \oplus U_1 \ar[r] & M_{2} \ar[r] & 0}.
		\end{equation*}
		So by  Corollary \ref{CorMono}, there is a monomorphism $M_1^{E(1)} \rightarrowtail M_2 \oplus U$. Since $M_1 < M_0$, our assumptions yield $M_2 \notin \mathcal{C}_{\underline{w}}$.

		By the same argument, we can see that $M_3, M_4, \dots, M_{m_{ij}-1} \notin \mathcal{C}_{\underline{w}}$ and by (S4) this means $m_{ij} < \infty$.
		
		By Lemma \ref{LemWordMij}, if we choose $\underline{u}$ and $\underline{v}$ so that $\underline{w} = \underline{u}s_i \underline{v}$ and there is some $\rho$ so that $\rho(\underline{w}) = \rho(\underline{u})(s, i)\rho$, then
		$\underline{w} \equiv \underline{u}' \{s_is_j\}^{m_{ij}} \underline{v}$ for some word $\underline{u}'$. 
		
		We still need to show that $\underline{u}' \{s_is_j\}^{m_{ij}} \underline{v}$ is equivalent to a word which is smaller than $\underline{w}$. 
		
		To do this, we use Lemma \ref{LemSame}. Either there is nothing to show, or there are $\rho_1, \dots, \rho_4$ and a pair $(r, h)$ so that
		\begin{equation} \label{Form5}
		\rho(\underline{w}) = \rho_1 (r, h) \rho_2 (s, i) \rho_3
		\end{equation}
		and there is some $\underline{w}'' \equiv \underline{w}$ with
		\begin{equation} \label{Form6}
		\rho(\underline{w}'') = \rho_1 \rho_2 (s, i) (t+1, j)\rho_3.
		\end{equation}
		We can assume that the word $\underline{x}$ with $\rho(\underline{x}) = \rho_1 (r, h)\rho_2$ is reduced, because otherwise there is nothing to prove. So by Lemma \ref{Lem2Ind}, there are some sequences of modules as in Lemma \ref{LemExSeq} that yield some $U' \in \mathcal{I}$ and an exact sequence
		\begin{equation*}
		\xymatrix{0 \ar[r] & M_0 \ar[r] & M_{m_{ij}-2} \oplus U' \ar[r] & Y \ar[r] &0}
		\end{equation*}
		with either $Y \in \mathcal{C}_{\underline{w}}$ or $Y = (r, h)^{E(m_{ij}-3)}$.
		
		By Corollary \ref{CorMono}, there is a monomorphism
		\begin{equation*}
		M_{m_{ij}-2} \rightarrowtail U \oplus Y \in \mathcal{C}_{\underline{w}''}.
		\end{equation*}
		By (\ref{Form5}), (\ref{Form6}) and the induction hypothesis, $\underline{w}''$ is not leftmost.
		
		So there is some $\underline{w}_3 \equiv \underline{w}''$ with $\underline{w}_3 <_l \underline{w}''$. We still need to show that $\underline{w}_3 <_l \underline{w}$.
		
		If $\underline{w}''$ is not reduced, this is obvious. If $Y \in \mathcal{C}_{\underline{w}}$, we can use the inductive assumption.
		
		So suppose that $\underline{w}''$ is reduced and $Y = (r, h)^{E(m_{ij}-3)}$. We denote the sequences of modules that fulfil (S1) - (S5) with respect to $M_0$ and $U$ by
		\begin{align*}
		\begin{split}
		&(X_1, X_2, \dots, X_{m})\\
		&(X'_1, X'_2,\dots, X'_{m})\\
		&(Y_1, Y_2, \dots, Y_{m}).
		\end{split}
		\end{align*}	
		There are sequences of modules
		\begin{align*}
		\begin{split}
		&(^1X_1, ^1X_2, \dots, ^1X_{m'})\\
		&(^1X'_1, ^1X'_2,\dots, ^1X'_{m'})\\
		&(^1Y_1, ^1Y_2, \dots, ^1Y_{m'})
		\end{split}
		\end{align*}
		that fulfil (S1) - (S5) with respect to $M_{m_{ij}-2}$ and $U \oplus (r, h)^{E(m_{ij}-3)}$. Let $(r', h')$ be the smallest indecomposable direct summand of $^1Y_{m'}$. Then $(r', h') < (r, h)$.
		
		If there is a pair $(r'', h'') \neq  (r', h')$ and series of pairs $\rho'_1, \rho'_2, \rho'_3$ so that we can write
		\begin{equation*}
		\rho(\underline{w}_3) = \rho'_1 (r', h') \rho'_2 \quad \text{and} \quad 
		\rho(\underline{w}) = \rho'_1 (r'', h'')\rho'_3,
		\end{equation*}
		then $\underline{w}_3 <_l \underline{w}''$ implies $\underline{w}_3 <_l \underline{w}$.
		
		A simple induction on $m'$ shows that this is indeed the case: If $m' = 1$, then $(r, h) = M_{m_{ij}-1}$, $(r', h') = M_{m_{ij}}$ and the assertion is true by Corollary \ref{CorStartInd}.
		
		By \ref{LemAlgInd}, the smallest direct summand of $^1Y_{m'}$ is also the smallest direct summand of $Y_m$ and the inductive step is obvious.
		
		So $\underline{w}$ is not leftmost and the proof is complete.
	\end{proof}
	
	\section{The other direction} \label{SecProof3}
	
	In this section we finally conclude the proof that the map between words and full additive cofinite subcategories of $\mod A$ introduced in Definition \ref{DefBijection} gives rise to a bijection between the leftmost words and the cofinite submodule closed subcategories. Since every element of the Weyl group has a unique leftmost element, this gives a bijection between the Weyl group elements and the cofinite, submodule closed subcategories. 
	
	Again, we start with the basis of an induction:

	\begin{l1} \label{LemBegin}
		\begin{enumerate}[(a)]
			\item Let $\underline{w} := \underline{u} \{s_is_j\}^{m_{ij}} \underline{v}$. If $\mathcal{C}_{\underline{w}}$ is submodule closed, then
			\begin{equation*}
			\underline{w} <_l \underline{u} \{s_js_i\}^{m_{ij}} \underline{v}.
			\end{equation*}
			
			\item The category $\mathcal{C}_{\underline{u}s_is_i\underline{v}}$ is not submodule closed.
			
			\item Let $\underline{w}' := \underline{u} \{s_js_i\}^{m_{ij}+1} \underline{v}$. Then $\mathcal{C}_{\underline{w}'}$ is not submodule closed. 
		\end{enumerate}
	\end{l1}
	
	\begin{proof}
		We prove (a) by contraposition. By Definition \ref{DefBijection}, $\ind \mathcal{I} \setminus \mathcal{C}_{\underline{w}}$ consists of the modules which correspond to the reflections in $\underline{w}$. 
		
		Assume that 
		\begin{equation*}
		\underline{u} \{s_js_i\}^{m_{ij}} \underline{v} <_l \underline{u} \{s_is_j\}^{m_{ij}} \underline{v}  = \underline{w}
		\end{equation*}	
		and
		\begin{equation*}
		\rho(\underline{w})	= \rho (\underline{u}) \underbrace{(p,i) (q,j) (p+1,i) \dots}_{m_{ij} \text{ pairs}} \rho_1
		\end{equation*}
		for a sequence of pairs $\rho_1$.
		
		By Lemma \ref{LemOneExchange}, the module $(q-1, j)$ exists and by Definition \ref{DefBijection}, $\mathcal{C}_{\underline{w}}$ contains all indecomposable, preinjective modules $M$ with $(q-1, j) < M <(p,i)$ or $(p,i) < M < (q, j)$. 
		
		First, suppose that $m_{ij} = 2$. In this case, $\mathcal{C}_{\underline{w}}$ contains the middle term of the AR-sequence that starts in $(q,j)$ by Lemmas \ref{LemARNumbers} and \ref{LemARWeyl}. Since $(q, j) \notin \mathcal{C}_{\underline{w}}$, the subcategory is not submodule closed.

		Now let $m_{ij} \ge 3$. In \ref{DefMs}, we defined $M_0 := (s, i)$, $M_1 := (t,j), \dots$ for some arbitrary, fixed $s, t$. By Remark \ref{RemRecSec}, we can assume without loss of generality that $m_{ij}$ is odd and we can choose $s, t$ so that $M_{m_{ij}-1} = (p, i)$.
		
		Then $M_{m_{ij}} = (q-1, j) \neq 0$ and by Lemma \ref{LemRec}, $E(m_{ij}-2) \neq 0$. By Lemma \ref{LemRecSeq}, there is an exact sequence
		\begin{equation*}
		\xymatrix{0 \ar[r] & M_0 \ar[r] & (M_{m_{ij}-1})^{E(m_{ij}-1)} \oplus U_{m_{ij}-1} \ar[r] & M_{m_{ij}}^{E(m_{ij}-2)} \ar[r] & 0}
		\end{equation*}
		so that no $M_0, M_1, \dots, M_{m_{ij}-1}$ is a direct summand of $U_{m_{ij}-1}$. By Lemma \ref{LemRec} we have $E(m_{ij}-1)=0$, so there is a monomorphism
		\begin{equation*}
		M_0 \rightarrowtail U_{m_{ij}-1}. 
		\end{equation*}
		
		It remains to show that $U_{m_{ij}-1} \in \mathcal{C}_{\underline{w}_2}$ by Corollary \ref{CorRecSeq}: If $X$ be a direct summand of $U_{m_{ij}-1}$, then $M_{m_{ij}} < X < M_0$ and thus $X \in \mathcal{C}_{\underline{w}}$. 
		
		By Lemma \ref{LemARNumbers}, part (b) is obvious.
		
		The proof of (c) is completely analogous to the proof of (a). 
	\end{proof}
	
	Finally, we are prepared to prove that Definition \ref{DefBijection} gives a bijection. Recall that $\mathcal{C}_w = \mathcal{C}_{\underline{w}}$ if $\underline{w}$ is the leftmost word for $w$:
	
	\begin{t1} \label{Th}
		The map $w \mapsto \mathcal{C}_w$ is a bijection between the elements of the Weyl group of $A$ and the cofinite submodule closed subcategories of $\mod A$.
	\end{t1}
	
	\begin{proof}
		The map is well defined by Lemma \ref{LemFirstDirection} and obviously injective. It remains to prove that it is surjective, i.e. that for all cofinite submodule closed subcategories $\mathcal{C}$ of $\mod A$, there is a $w \in W$, so that $\mathcal{C} = \mathcal{C}_w$. 
		
		We already know that $\mathcal{C} = \mathcal{C}_{\underline{w}}$ for some word from Lemma \ref{LemWord}, so we only need to show that $\underline{w}$ is leftmost. 
		
		Assume that the word $\underline{w}$ for the element $w \in W$ is not leftmost. We show that $\mathcal{C}_{\underline{w}}$ is not submodule closed by induction on the number of Coxeter relations that are needed to transform $\underline{w}$ into a smaller word. 
		
		If only one relation is needed, then the theorem is the result of Lemma \ref{LemBegin}. Now suppose that the assertion is true if we need $m$ or less relations and that we need $m+1$ relations to transform $\underline{w}$ into a smaller word.
		
		Then there are some $1 \le i, j \le n$ and some words $\underline{x}$, $\underline{x}'$, $\underline{y}$ so that
		\begin{equation} \label{Form4}
		\underline{w} = \underline{x} s_i \underline{y} \equiv \underline{x}' s_j \underline{y} = \underline{w}'
		\end{equation}
		and
		\begin{equation*}
		\underline{w}' <_l \underline{w}
		\end{equation*}
		with $i \neq j$.
		
		Thus, there are some words $\underline{w}'', \underline{x}''$ so that
		\begin{equation*}
		\underline{w} \equiv \underline{w}'' = \underline{x}'' \{s_is_j\}^{m_{ij}} \underline{y}.
		\end{equation*}
		Because of the inductive assumption, we can suppose that $\underline{w} \le_l \underline{w}''$ and that $\underline{x}$ is leftmost. Obviously, we can choose $\underline{x}''$ to be leftmost.
		
		Let the reflection $s_i$ in (\ref{Form4}) correspond to $M_0$. We can assume that $m_{ij} \ge 3$, since there is nothing to show if the middle term of the Auslander-Reiten sequence that starts in $M_0$ is contained in $\mathcal{C}_{\underline{w}}$. We can also assume that there is some word $\underline{x}''$ so that $\underline{x} = \underline{x}''s_j$: If there are $s_{k_1}, \dots s_{k_m}$ with $m_{k_1, i} = \dots = m_{k_m, i}$ and $\underline{x} = \underline{x}''s_js_{k_1} \dots s_{k_m}$, then $\underline{x}s_i \equiv  \underline{x}'' s_js_is_{k_1} \dots s_{k_m}$ and if there is some $U \in \mathcal{C}_{\underline{x}''s_js_i}$ with a monomorphism $M_0 \rightarrowtail U$, then $U \in \mathcal{C}_{\underline{w}}$.
		
		Without loss of generality, we can assume that $m_{ij}$ is odd; otherwise we relabel $i$ and $j$ and get the same arguments by Remark \ref{RemRecSec}.
		
		By Lemma \ref{LemPart2}, we can suppose that $M_0, M_1, \dots, M_{m_{ij}-3}, M_{m_{ij}-2} \in \mathcal{C}_{\underline{w}}$. We show that there is some $U \in \mathcal{C}_{\underline{w}}$ with a monomorphism $M_0 \rightarrowtail U$.
		
		Let $(q, j):= M_{m_{ij}-2}$. We use Lemma \ref{LemSame}: Because $m > 1$, there is some $\underline{w}_3 \equiv \underline{w}$ with series of pairs $\rho_1, \dots, \rho_4$ and a pair $(r, h)$ so that
		\begin{equation} \label{Word1}
		\rho(\underline{w}) = \rho_1 (r, h) \rho_2 (q, j) \rho_3 (s, i) \rho_4
		\end{equation}
		and
		\begin{equation} \label{Word2}
		\rho(\underline{w}_3) = \rho_1 \rho_2 (q, j) \rho_3 (s, i) (t+1, j)\rho_4.
		\end{equation}

		By Lemma \ref{Lem2Ind}, if $m_{ij} \ge 3$, then there is an exact sequence	
		\begin{equation*}
		\xymatrix{0 \ar[r] &M_0 \ar[r] &M_{m_{ij}-2}\oplus U'\ar[r] & Y \ar[r] &0}
		\end{equation*}
		so that either $Y \in \mathcal{C}_{\underline{w}}$ or both $Y = (r, h)^{E(m_{ij}-3)}$ and $U' \in \mathcal{C}_{\underline{w}}$ hold.
	
		We want to show that there is some $U'' \in \mathcal{C}_{\underline{w}}$ and a monomorphism
		\begin{equation} \label{Equation2}
		M_{m_{ij}-2} = (q, j) \rightarrowtail U'' \oplus Y.
		\end{equation}
		We prove this inductively: First, note that the word $\underline{x}'' s_j$ is not leftmost Lemma \ref{LemIndStep} and \ref{RemIndStep}.
		
		If $\underline{w} = \underline{w}''$, then $(r, h) = M_{m_{ij}-1}$. So by the inductive hypothesis, there is a monomorphism $(q, j) \rightarrowtail U'' \oplus (r, h)^{\gamma}$.
		
		Since $E(1) = E(m_{ij}-3)$ by table (\ref{Table}), we obviously get $Y = (r, h)^{E(m_{ij}-3)}$ and $\gamma = E(m_{ij}-3)$. 
		
		The inductive step is completely analogous to the one in Lemma \ref{Lem2Ind}.
		
		By our assumptions, $\underline{x}$ is leftmost and thus $Y \notin \mathcal{C}_{\underline{w}}$ by Lemma \ref{LemFirstDirection}. So $U' \in \mathcal{C}_{\underline{w}}$.
		
		By Lemma \ref{LemAlgInd}, there is a monomorphism $M_0 \rightarrowtail U' \oplus U'' \in \mathcal{C}_{\underline{w}}$ and $\mathcal{C}_{\underline{w}}$ is not submodule closed.
	\end{proof}

	\section{Some consequences} \label{SecConsequences}
	
	We conclude the chapter with a generalization and a corollary:
	
	As in \cite{ORT}, Section 8 we can extend the notion of leftmost words:
	
	\begin{d1} \label{DefInf}
		Define infinite words analogously to words, only as infinite instead of finite sequences. We say that an (infinite) word is \textit{leftmost} if any initial subword of finite length is leftmost.
	\end{d1}
	
	Analogously to \cite{ORT}, Theorem 8.1, we get the following:
	
	\begin{t1}
		There is a bijection between the (finite and infinite) leftmost words over $S = \{s_1, s_2, \dots s_n\}$ and the submodule closed subcategories of $\mathcal{I}$, the preinjective component of $\mod A$.
	\end{t1}
	
	\begin{proof}
		This is completely analogous to \cite{ORT}, 8.1: 
		
		Let $\mathcal{C}$ be a submodule closed subcategory of $\mathcal{I}$. Since $\mathcal{I}$ contains at most a countable number of indecomposable modules, we can set 
		\begin{equation*}
		\ind \mathcal{I} \setminus \mathcal{C} =: \{(r_1, i_1), (r_2, i_2), \dots\}
		\end{equation*}
		and $(r_1, i_1) < (r_2, i_2) < \dots$. Then the subcategory $\mathcal{C}_m$ with
		\begin{equation*}
		\ind \mathcal{I} \setminus \mathcal{C}_m = \{(r_1, i_1), (r_2, i_2), \dots, (r_m, i_m)\}
		\end{equation*}
		is submodule closed for all $m \in \mathbb{N}$. By Lemma \ref{LemPreinj} and Theorem \ref{Th}, the words $\underline{w}_m$ with 
		\begin{equation*}
		\rho(\underline{w}_m) = (r_1, i_1) (r_2, i_2) \dots (r_m, i_m)
		\end{equation*}
		are leftmost for all $m \in \mathbb{N}$. By Definition \ref{DefInf}, the (infinite) word $\underline{w}$ with 
		\begin{equation*}
		\rho(\underline{w}) = (r_1, i_1)(r_2, i_2) \dots
		\end{equation*}
		is leftmost and $\mathcal{C} = \mathcal{C}_{\underline{w}}.$
		
		On the other hand assume that the (infinite) word $\underline{w}$ with 
		\begin{equation*}
		\rho(\underline{w}) = (r_1, i_1)(r_2, i_2) \dots
		\end{equation*}
		is leftmost. Then 
		the words with
		\begin{equation*}
		\rho(\underline{w}_m) = (r_1, i_1) (r_2, i_2) \dots (r_m, i_m)
		\end{equation*}
		are leftmost for all $m \in \mathbb{N}$. By \ref{Th}, the categories $\mathcal{C}_{\underline{w}_m}$ are submodule closed. Thus $\mathcal{C}_{\underline{w}}$ is also submodule closed: if there was a module $M \notin \mathcal{C}_{\underline{w}}$ and some module $U \in \mathcal{C}_{\underline{w}}$ with a monomorphism $M \rightarrowtail U$, then $U \in \mathcal{C}_{\underline{w}_m}$ for all $m \in \mathbb{N}$ and there is some $m \in \mathbb{N}$ so that $M \notin \mathcal{C}_{\underline{w}_m}$, since $M$ is finitely generated.
	\end{proof}
	
	We can draw a further corollary. Let $A'$ be a hereditary and let the module category $\mod A'$ be equivalent to the subcategory of $\mod A$ with the simple modules $S_j$, $j \in J$ for some $J \subseteq \{1, 2, \dots n\}$. Let $\mathcal{I}_{A'}$ be the subcategory of $\mod A'$ consisting of all preinjective modules.
	
	\begin{c1}
		There is a bijection between the submodule closed subcategories of $\mathcal{I}_{A'}$ and the submodule closed subcategories of $\mathcal{I}$ which contain all $\tau^r I_i$ with $r\in \mathbb{N}_0$ and $i \in \mathbb{N} \setminus J$.
	\end{c1}
	
	\begin{proof}
		The words in the Weyl group of $A'$ are exactly the words in the Weyl group of $A$ which only consist of reflections $s_j$ with $j \in J$.
	\end{proof}

\section*{Acknowledgements}

The author would like to thank Henning Krause for his suggestions and insights and Hugh Thomas for his comments concerning the exposition.

\end{document}